\documentclass{gtart}
\usepackage{amsfonts}
\usepackage{amsmath}

\input gtoutput
\volumenumber{4}\papernumber{6}\volumeyear{2000}
\pagenumbers{179}{218}
\proposed{Jean-Pierre Otal}
\seconded{Walter Neumann, Joan Birman}
\received{18 May 1999}
\accepted{24 July 2000}
\revised{6 April 2000}
\published{9 August 2000}

\newtheorem{theorem}{Theorem}[section]
\newtheorem{conjecture}[theorem]{Conjecture}
\newtheorem{corollary}[theorem]{Corollary}
\newtheorem{lemma}[theorem]{Lemma}
\newtheorem{proposition}[theorem]{Proposition}

\theoremstyle{remark}
\newtheorem{definition}[theorem]{Definition}

\newtheorem{remark}[theorem]{Remark}

\begin{document}

\title{Splittings of groups and intersection numbers}

\author{Peter Scott\\Gadde A Swarup}

\address{Mathematics Department, University of Michigan\\
Ann Arbor, Michigan 48109, USA\\{\rm and}\\
Mathematics Department, University of Melbourne\\
Parkville, Victoria 3052, Australia}

\asciiaddress{Mathematics Department, University of Michigan\\
Ann Arbor, Michigan 48109, USA\\
Mathematics Department, University of Melbourne\\
Parkville, Victoria 3052, Australia}

\email{pscott@math.lsa.umich.edu\\gadde@ms.unimelb.edu.au}

\begin{abstract}
We prove algebraic analogues of the facts that a curve on a surface with
self-intersection number zero is homotopic to a cover of a simple curve, and
that two simple curves on a surface with intersection number zero can be
isotoped to be disjoint.
\end{abstract}




\maketitlepage

In this paper, we will discuss an algebraic version of intersection numbers
which was introduced by Scott in \cite{Scott:Intersectionnumbers}. First we
need to discuss intersection numbers in the topological setting. Let $F$
denote a surface and let $L$ and $S$ each be a properly immersed two-sided
circle or compact arc in $F$. Here `properly' means that the boundary of the 
$1$--manifold lies in the boundary of $F$. One can define the intersection
number of $L$ and $S$ to be the least number of intersection points
obtainable by homotoping $L$ and $S$ transverse to each other. (The count is
to be made without any signs attached to the intersection points.) It is
obvious that this number is symmetric in the sense that it is independent of
the order of $L$ and $S$. It is also obvious that $L$ and $S$ have
intersection number zero if and only if they can be properly homotoped to be
disjoint. It seems natural to define the self-intersection number of an
immersed two-sided circle or arc $L$ in $F$ to be the least number of
transverse intersection points obtainable by homotoping $L$ into general
position. With this definition, $L$ has self-intersection number zero if and
only if it is homotopic to an embedding. However, in light of later
generalisations, it turns out that this definition should be modified a
little in order to ensure that the self-intersection number of any cover of
a simple closed curve is also zero. No modification is needed unless $L$ is
a circle which can be homotoped to cover another immersion with degree
greater than $1$. In this case, suppose that the maximal degree of covering
which can occur is $k$ and that $L$ covers $L^{\prime }$ with degree $k$.
Then we define the self-intersection number of $L$ to be $k^{2}$ times the
self-intersection number of $L^{\prime }$. With this modified definition, $%
L\;$has self-intersection number zero if and only if it can be homotoped to
cover an embedding.

In \cite{FHSleastareasurfaces}, Freedman, Hass and Scott introduced a notion
of intersection number and self-intersection number for two-sided $\pi _{1}$%
--injective immersions of compact surfaces into $3$--manifolds which
generalises the preceding ideas. Their intersection number cannot be
described as simply as for curves on a surface, but it does share some
important properties. In\ particular, it is a non-negative integer and it is
symmetric, although this symmetry is not obvious from the definition.
Further, two surfaces have intersection number zero if and only if they can
be homotoped to be disjoint, and a single surface has self-intersection
number zero if and only if it can be homotoped to cover an embedding. These
two facts are no longer obvious consequences of the definition, but are
non-trivial applications of the theory of least area surfaces.

In \cite{Scott:Intersectionnumbers}, Scott extended the ideas of \cite
{FHSleastareasurfaces} to define intersection numbers in a purely group
theoretic setting. The details will be discussed in the first section of
this paper, but we give an introduction to the ideas here. It seems clear
that everything discussed in the preceding two paragraphs should have a
purely algebraic interpretation in terms of fundamental groups of surfaces
and $3$--manifolds, and the aim is to find an interpretation which makes
sense for any group. It seems natural to attempt to define the intersection
number of two subgroups $H\;$and $K$ of a given group $G$. This is exactly
what the topological intersection number of simple closed curves on a
surface does when $G$ is the fundamental group of a closed orientable
surface and we restrict attention to infinite cyclic subgroups $H$ and $K$.
However, if one considers two simple arcs on a surface $F$ with boundary,
they each carry the trivial subgroup of $G=\pi _{1}(F)$, whereas we know
that some arcs have intersection number zero and others do not. Thus
intersection numbers are not determined simply by the groups involved. We
need to look a little deeper in order to formulate the algebraic analogue.
First we need to think a bit more about curves on surfaces. Let $L$ be a
simple arc or closed curve on an orientable surface $F$, let $G$ denote $\pi
_{1}(F)$ and let $H$ denote the image of $\pi _{1}(L)$ in $G$. If $L$
separates $F$ then, in most cases, it gives $G$ the structure of an
amalgamated free product $A\ast _{H}B$, and if $L$ is non-separating, it
gives $G$ the structure of a HNN extension $A\ast _{H}$. In order to avoid
discussing which of these two structures $G$ has, it is convenient to say
that a group $G$ \textit{splits over a subgroup} $H$ if $G$ is isomorphic to 
$A\ast _{H}$ or to $A\ast _{H}B$, with $A\neq H\neq B$. (Note that the
condition that $A\neq H\neq B$ is needed as otherwise any group $G$ would
split over any subgroup $H$. For one can always write $G=G\ast _{H}H$.)
Thus, in most cases, $L$ determines a splitting of $G=\pi _{1}(F)$. Usually
one ignores base points, so that the splitting of $G$ is only determined up
to conjugacy. In \cite{Scott:Intersectionnumbers}, Scott defined the
intersection number of two splittings of any group $G$ over any subgroups $H$
and $K$. In the special case when $G$ is the fundamental group of a compact
surface $F$ and these splittings arise from embedded arcs or circles on $F$,
the algebraic intersection number of the splittings equals the topological
intersection number of the corresponding $1$--manifolds. The analogous
statement holds when $G$ is the fundamental group of a compact $3$--manifold
and these splittings arise from $\pi _{1}$--injective embedded surfaces. In
general, the algebraic intersection number shares some properties of the
topological intersection number. Algebraic intersection numbers are
symmetric, and if $G$, $H$ and $K$ are finitely generated, the intersection
number of splittings of $G$ over $H$ and over $K$ is a non-negative integer.

The first main result of this paper is a generalisation to the algebraic
setting of the fact that two simple arcs or closed curves on a surface have
intersection number zero if and only if they can be isotoped apart. Of
course, the idea of isotopy makes no sense in the algebraic setting, so we
need some algebraic language to describe multiple disjoint curves on a
surface. Let $L_{1},\ldots ,L_{n}$ be disjoint simple arcs or closed curves
on a compact orientable surface $F$ with fundamental group $G$, such that
each $L_{i}$ determines a splitting of $G$. Together they determine a graph
of groups structure on $G$ with $n$ edges. We say that a collection of $n$
splittings of a group $G$ is \textit{compatible} if $G$ can be expressed as
the fundamental group of a graph of groups with $n$ edges, such that, for
each $i$, collapsing all edges but the $i$-th yields the $i$-th splitting of 
$G.$ We will say that the splittings are \textit{compatible up to conjugacy}
if collapsing all edges but the $i$-th yields a splitting of $G$ which is
conjugate to the $i$-th given splitting. Clearly disjoint essential simple
arcs or closed curves on $F$ define splittings of $G$ which are compatible
up to conjugacy. The precise statement we obtain is the following.

\ppar

\textbf{Theorem \ref{disjointsplittings}}\qua{\sl{Let }$G$\sl{\ be a
finitely generated group with }$n$\sl{\ splittings over finitely
generated subgroups. This collection of splittings is compatible up to
conjugacy if and only if each pair of splittings has intersection number
zero. Further, in this situation, the graph of groups structure on }$G$%
\sl{\ obtained from these splittings has a unique underlying graph, and
the edge and vertex groups are unique up to conjugacy.}}

\ppar

So far, we have not discussed any algebraic analogue of non-embedded arcs or
circles on surfaces. There is such an analogue which is the idea of an
almost invariant subset of the quotient $H\backslash G$, where $H$ is a
subgroup of $G$. This generalises the idea of an immersed curve in a surface
or of an immersed $\pi _{1}$--injective surface in a 3--manifold which carries
the subgroup $H$ of $G$. We give the definitions in section 1. There is also
an idea of intersection number of such things, which we give in Definition 
\ref{defnofintersectionnumber}. This too was introduced by Scott in \cite
{Scott:Intersectionnumbers}. Our second main result, Theorem \ref
{splittingsexist}, is an algebraic analogue of the fact that a singular
curve on a surface or a singular surface in a $3$--manifold which has
self-intersection number zero can be homotoped to cover an embedding. It
asserts that if $H\backslash G$ has an almost invariant subset with
self-intersection number zero, then $G$ has a splitting over a subgroup $%
H^{\prime }$ commensurable with $H$. We leave the precise statement until
section 2.

In a separate paper \cite{S-S2}, we use the ideas about intersection numbers
of splittings developed in \cite{Scott:Intersectionnumbers} and in this
paper to study JSJ decompositions of Haken $3$--manifolds. The problem there
is to recognize which splittings of the fundamental group of such a manifold
arise from the JSJ decomposition (see \cite{LS:Non-positive} and \cite
{NS:Canonical}). It turns out that a class of splittings which we call
canonical can be defined using intersection numbers and we use this to show
that the JSJ decomposition for Haken $3$--manifolds depends only on the
fundamental group. This leads to an algebraic proof of Johannson'
Deformation Theorem. It seems very likely that similar ideas apply to Sela's
JSJ decompositions \cite{Sela:JSJ} of hyperbolic groups and thus provide a
common thread to the two types of JSJ decomposition. Thus, the use of
intersection numbers seems to provide a tool in the study of diverse topics
in group theory and this paper together with \cite{Scott:Intersectionnumbers}
provides some of the foundational material.

This paper is organised as follows. In section 1, we recall from \cite
{Scott:Intersectionnumbers} the basic definitions of intersection numbers in
the algebraic context. We also prove a technical result which was
essentially proved by Scott \cite{Scott:TorusTheorem} in 1980. However,
Scott's results were all formulated in the context of surfaces in $3$%
--manifolds, so we give a complete proof of the generalisation to the purely
group theoretic context. Section 2 is devoted to the proofs of our two main
results discussed above.

There is a second natural idea of intersection number, which we discuss in
section 3. We call it the strong intersection number. It is not symmetric in
general, but this is not a problem when one is considering self-intersection
numbers. We also discuss when the two kinds of intersection number are
equal, which then forces the strong intersection number to be symmetric. We
use these ideas to give a new approach to a result of Kropholler and Roller 
\cite{KR:1} on splittings of Poincar\'{e} duality groups. We also discuss
applications of our ideas to prove a special case of a conjecture of
Kropholler and Roller \cite{KR:2} on splittings of groups in general. We
point out that these ideas lead to an alternative approach to the algebraic
Torus Theorem \cite{D-S:torus}. We end the section with a brief discussion
of an error in \cite{Scott:Intersectionnumbers}. In section 3 of that paper,
Scott gave an incorrect interpretation of the intersection number of two
splittings. His error was caused by confusing the ideas of strong and
ordinary intersection. However, the arguments in \cite
{Scott:Intersectionnumbers} work to give a nice interpretation of the
intersection number in the case when it is equal to the strong intersection
number. Without this condition, finding nice interpretations of the two
intersection numbers is an open problem.

\section{Preliminaries and statements of main results}

We will start by recalling from \cite{Scott:Intersectionnumbers} how to
define intersection numbers in the algebraic setting. We will connect this
with the natural topological idea of intersection number already discussed
in the introduction. Consider two simple closed curves $L$ and $S$ on a
closed orientable surface $F$. As in \cite{FHScurvesonsurfaces}, it will be
convenient to assume that $L$ and $S$ are shortest geodesics in some
Riemannian metric on $F$ so that they automatically intersect minimally. We
will interpret the intersection number of $L$ and $S$ in suitable covers of $%
F$, exactly as in \cite{FHScurvesonsurfaces} and \cite{FHSleastareasurfaces}%
. Let $G$ denote $\pi _{1}(F)$, let $H$ denote the infinite cyclic subgroup
of $G$ carried by $L$, and let $F_{H}$ denote the cover of $F$ with
fundamental group equal to $H$. Then $L$ lifts to $F_{H}$ and we denote its
lift by $L$ again. Let $l$ denote the pre-image of this lift in the
universal cover $\widetilde{F}$ of $F$. The full pre-image of $L$ in $%
\widetilde{F}$ consists of disjoint lines which we call $L$--lines, which are
all translates of $l$ by the action of $G$. (Note that in this paper groups
act on the left on covering spaces.) Similarly, we define $K$, $F_{K}$, the
line $s$ and $S$--lines in $\widetilde{F}$. Now we consider the images of the 
$L$--lines in $F_{K}$. Each $L$--line has image in $F_{K}$ which is a line or
circle. Then we define $d(L,S)$ to be the number of images of $L$--lines in $%
F_{K}$ which meet $S$. Similarly, we define $d(S,L)$ to be the number of
images of $S$--lines in $F_{H}$ which meet $L$. It is shown in \cite
{FHScurvesonsurfaces}, using the assumption that $L$ and $S$ are shortest
closed geodesics, that each $L$--line in $F_{K}$ crosses $S$ at most once,
and similarly for $S$--lines in $F_{H}$. It follows that $d(L,S)$ and $d(S,L)$
are each equal to the number of points of $L\cap S$, and so they are equal
to each other.

We need to take one further step in abstracting the idea of intersection
number. As the stabiliser of $l$ is $H$, the $L$--lines naturally correspond
to the cosets $gH$ of $H$ in $G$. Hence the images of the $L$--lines in $%
F_{K} $ naturally correspond to the double cosets $KgH$. Thus we can think
of $d(L,S)$ as the number of double cosets $KgH$ such that $gl$ crosses $s$.
This is the idea which we generalise to define intersection numbers in a
purely algebraic setting.

First we need some terminology.

Two sets $P$ and $Q$ are \textit{almost equal} if their symmetric difference 
$P-Q\cup Q-P$ is finite. We write $P\overset{a}{=}Q.$

If a group $G$ acts on the right on a set $Z$, a subset $P$ of $Z$ is 
\textit{almost invariant} if $Pg\overset{a}{=}P$ for all $g$ in $G$. An
almost invariant subset $P$ of $Z$ is \textit{non-trivial} if $P$ and its
complement $Z-P$ are both infinite. The complement $Z-P$ will be denoted
simply by $P^{\ast }$, when $Z$ is clear from the context

For finitely generated groups, these ideas are closely connected with the
theory of ends of groups via the Cayley graph $\Gamma $ of $G$ with respect
to some finite generating set of $G$. (Note that $G$ acts on its Cayley
graph on the left.) Using $\mathbb{Z}_{2}$ as coefficients, we can identify $%
0$--cochains and $1$--cochains on $\Gamma $ with sets of vertices or edges. A
subset $P$ of $G$ represents a set of vertices of $\Gamma $ which we also
denote by $P$, and it is a beautiful fact, due to Cohen \cite{Cohen}, that $%
P $ is an almost invariant subset of $G$ if and only if $\delta P$ is
finite, where $\delta $ is the coboundary operator. Now $\Gamma $ has more
than one end if and only if there is an infinite subset $P$ of $G$ such that 
$\delta P $ is finite and $P^{\ast }$ is also infinite. Thus $\Gamma $ has
more than one end if and only if $G$ contains a non-trivial almost invariant
subset. If $H$ is a subgroup of $G$, we let $H\backslash G$ denote the set
of cosets $Hg$ of $H$ in $G$, ie, the quotient of $G$ by the left action of 
$H$. Of course, $G$ will no longer act on the left on this quotient, but it
will still act on the right. Thus we also have the idea of an almost
invariant subset of $H\backslash G$, and the graph $H\backslash \Gamma $ has
more than one end if and only if $H\backslash G$ contains a non-trivial
almost invariant subset. Now the number of ends $e(G)$ of $G$ is equal to
the number of ends of $\Gamma $, so it follows that $e(G)>1$ if and only if $%
G$ contains a non-trivial almost invariant subset. Similarly, the number of
ends $e(G,H)$ of the pair $(G,H)$ equals the number of ends of $H\backslash
\Gamma $, so that $e(G,H)>1$ if and only if $H\backslash G$ contains a
non-trivial almost invariant subset.

Now we return to the simple closed curves $L$ and $S$ on the surface $F$.
Pick a generating set for $G$ which can be represented by a bouquet of
circles embedded in $F$. We will assume that the wedge point of the bouquet
does not lie on $L$ or $S$. The pre-image of this bouquet in $\widetilde{F}$
will be a copy of the Cayley graph $\Gamma $ of $G$ with respect to the
chosen generating set. The pre-image in $F_{H}$ of the bouquet will be a
copy of the graph $H\backslash \Gamma $, the quotient of $\Gamma $ by the
action of $H$ on the left. Consider the closed curve $L$ on $F_{H}$. Let $P$
denote the set of all vertices of $H\backslash \Gamma $ which lie on one
side of $L$. Then $P$ has finite coboundary, as $\delta P$ equals exactly
the edges of $H\backslash \Gamma $ which cross $L$. Hence $P$ is an almost
invariant subset of $H\backslash G$. Let $X$ denote the pre-image of $P$ in $%
\Gamma $, so that $X$ equals the set of vertices of $\Gamma $ which lie on
one side of the line $l$. Now finally the connection between the earlier
arguments and almost invariant sets can be given. For we can decide whether
the lines $l$ and $s$ cross by considering instead the sets $X$ and $Y$. The
lines $l$ and $s$ together divide $G$ into the four sets $X\cap Y$, $X^{\ast
}\cap Y$, $X\cap Y^{\ast }$ and $X^{\ast }\cap Y^{\ast }$, where $X{\ }%
^{\ast }$ denotes $G-X$, and $l$ crosses $s$ if and only if each of these
four sets projects to an infinite subset of $K\backslash G.$

Now let $G$ be a group with subgroups $H$ and $K$, let $P$ be a non-trivial
almost invariant subset of $H\backslash G$ and let $Q$ be a non-trivial
almost invariant subset of $K\backslash G$. We will define the intersection
number $i(P,Q)$ of $P$ and $Q$. First we need to consider the analogues of
the sets $X\;$and $Y$ in the preceding paragraph, and to say what it means
for them to cross.

\begin{definition}
If $G$ is a group and $H$ is a subgroup, then a subset $X$ of $G$ is $H$%
\textsl{-almost invariant} if $X$ is invariant under the left action of $H$,
and simultaneously $H\backslash X$ is an almost invariant subset of $%
H\backslash G$. In addition, $X$ is a \textsl{non-trivial} $H$--almost
invariant subset of $G$, if the quotient sets $H\backslash X$ and $%
H\backslash X^{\ast }$ are both infinite.
\end{definition}

Note that if $H$ is trivial, then a $H$--almost invariant subset of $%
G$ is the same as an almost invariant subset of $G$.

\begin{definition}
Let $X$ be a $H$--almost invariant subset of $G$ and let $Y$ be a $K$--almost
invariant subset of $G$. We will say that $X$ \textsl{crosses} $Y$ if each
of the four sets $X\cap Y$, $X^{\ast }\cap Y$, $X\cap Y^{\ast }$ and $%
X^{\ast }\cap Y^{\ast }$ projects to an infinite subset of $K\backslash G.$
\end{definition}

We will often write $X^{(\ast )}\cap Y^{(\ast )}$ instead of listing the
four sets $X\cap Y$, $X^{\ast }\cap Y$, $X\cap Y^{\ast }$ and $X^{\ast }\cap
Y^{\ast }.$

If $G$ is a group and $H$ is a subgroup, then we will say that a subset $W$
of $G$ is $H$--\textit{finite} if it is contained in the union of
finitely many left cosets $Hg$ of $H$ in $G$, and we will say that two
subsets $V$ and $W$ of $G$ are $H$--\textit{almost equal} if their
symmetric difference is $H$--finite.

In this language, $X$ crosses $Y$ if each of the four sets $X^{(\ast )}\cap
Y^{(\ast )}$ is not $K$--finite.

This definition of crossing is not symmetric, but it is shown in \cite
{Scott:Intersectionnumbers} that if $G$ is a finitely generated group with
subgroups $H$ and $K$, and $X$ is a non-trivial $H$--almost invariant subset
of $G$ and $Y$ is a non-trivial $K$--almost invariant subset of $G$, then $X$
crosses $Y$ if and only if $Y$ crosses $X$. If $X$ and $Y$ are both trivial,
then neither can cross the other, so the above symmetry result is clear.
However, this symmetry result fails if only one of $X$ or $Y$ is trivial.
This lack of symmetry will not concern us as we will only be interested in
non-trivial almost invariant sets.

Now we come to the definition of the intersection number of two almost
invariant sets.

\begin{definition}
\label{defnofintersectionnumber}Let $H$ and $K$ be subgroups of a finitely
generated group $G$. Let $P$ denote a non-trivial almost invariant subset of 
$H\backslash G$, let $Q$ denote a non-trivial almost invariant subset of $%
K\backslash G$ and let $X$ and $Y$ denote the pre-images of $P$ and $Q$
respectively in $G$. Then the intersection number $i(P,Q)$ of $P$ and $Q$
equals the number of double cosets $KgH$ such that $gX$ crosses $Y.$
\end{definition}

\begin{remark}
\label{almostequalsetshavesameintersectionnumber}The following facts about
the intersection number are proved in \cite{Scott:Intersectionnumbers}.

\begin{enumerate}
\item  Intersection numbers are symmetric, ie $i(P,Q)=i(Q,P)$.

\item  $i(P,Q)$ is finite when $G$, $H$, and $K$ are all finitely generated.

\item  If $P^{\prime }$ is an almost invariant subset of $H\backslash G$
which is almost equal to $P$ or to $P^{\ast }$ and if $Q^{\prime }$ is an
almost invariant subset of $K\backslash G$ which is almost equal to $Q$ or
to $Q^{\ast }$, then $i(P^{\prime },Q^{\prime })=i(P,Q).$
\end{enumerate}
\end{remark}

We will often be interested in situations where $X$ and $Y$ do not cross
each other and neither do many of their translates. This means that one of
the four sets $X^{(\ast )}\cap Y^{(\ast )}$ is $K$--finite, and similar
statements hold for many translates of $X\;$and $Y$. If $U=uX$ and $V=vY$ do
not cross, then one of the four sets $U^{(\ast )}\cap V^{(\ast )}$ is $K^{v}$%
--finite, but probably not $K$--finite. Thus one needs to keep track of which
translates of $X$ and $Y$ are being considered in order to have the correct
conjugate of $K$, when formulating the condition that $U$ and $V$ do not
cross. The following definition will be extremely convenient because it
avoids this problem, thus greatly simplifying the discussion at certain
points.

\begin{definition}
Let $U$ be a $H$--almost invariant subset of $G$ and let $V$ be a $K$--almost
invariant subset of $G$. We will say that $U\cap V$ is \textsl{small} if it
is $H$--finite.
\end{definition}

\begin{remark}
As the terminology is not symmetric in $U$ and $V$ and makes no reference to 
$H$ or $K$, some justification is required. If $U$ is also $H^{\prime }$%
--almost invariant for a subgroup $H^{\prime }$ of $G$, then $H^{\prime }$
must be commensurable with $H$. Thus $U\cap V$ is $H$--finite if and only if
it is $H^{\prime }$--finite. In addition, the fact that crossing is symmetric
tells us that $U\cap V$ is $H$--finite if and only if it is $K$--finite. This
provides the needed justification of our terminology.

Finally, the reader should be warned that this use of the word small has
nothing to do with the term small group which means a group with no
subgroups which are free of rank $2$.
\end{remark}

At this point we have the machinery needed to define the intersection number
of two splittings. This definition depends on the fact, which we recall from 
\cite{Scott:Intersectionnumbers}, that if a group $G$ has a splitting over a
subgroup $H$, there is a $H$--almost invariant subset $X$ of $G$ associated
to the splitting in a natural way. This is entirely clear from the
topological point of view as follows. If $G=A\ast _{H}B$, let $N$ denote a
space with fundamental group $G$ constructed in the usual way as the union
of $N_{A}$, $N_{B}$ and $N_{H}\times I$. If $G=A\ast _{H}$, then $N$ is
constructed from $N_{A}$ and $N_{H}\times I$ only. Now let $M$ denote the
based cover of $N$ with fundamental group $H$, and denote the based lift of $%
N_{H}\times I$ into $M$ by $N_{H}\times I$. Then $X$ corresponds to choosing
one side of $N_{H}\times I$ in $M$. We now give a purely algebraic
description of this choice of $X$ (see \cite{Scott-Wall:Topological} for
example). If $G=A\ast _{H}B$, choose right transversals $T_{A}$, $T_{B}$ of $%
H$ in $A$, $B$, both of which contain the identity element. (A right
transversal for a subgroup $H$ of a group $G$ consists of one representative
element for each right coset $gH$ of $H$ in $G.)$ Each element of $G$ can be
expressed uniquely in the form $a_{1}b_{1}a_{2}...a_{n}b_{n}h$ with $h\in H$%
, $a_{i}\in T_{A}$, $b_{i}\in T_{B}$, where only $h$, $a_{1}$ and $b_{n}$
are allowed to be trivial. Then $X$ consists of elements for which $a_{1}$
is non-trivial. In the case of a HNN--extension $A\ast _{H}$, let $\alpha
_{i} $, $i=1$, $2$, denote the two inclusions of $H$ in $A$ so that $%
t^{-1}\alpha _{1}(h)t=\alpha _{2}(h)$, and choose right transversals $T_{i}$
of $\alpha _{i}(H)$ in $A$, both of which contain the identity element. Each
element of $G$ can be expressed uniquely in the form $a_{1}t^{\epsilon
_{1}}a_{2}t^{\epsilon _{2}}...a_{n}t^{\epsilon _{n}}a_{n+1}$ where $a_{n+1}$
lies in $A$ and, for $1\leq i\leq n$, $\epsilon _{i}=1$ or $-1$, $a_{i}\in
T_{1}$ if $\epsilon _{i}=1$, $a_{i}\in T_{2}$ if $\epsilon _{i}=-1$ and
moreover $a_{i}\neq 1$ if $\epsilon _{i-1}\neq \epsilon _{i}$. In this case, 
$X$ consists of elements for which $a_{1}$ is trivial and $\epsilon _{1}=1$.
In both cases, the stabiliser of $X$ under the left action of $G$ is exactly 
$H$ and, for every $g\in G$, at least one of the four sets $X^{(\ast )}\cap
gX^{(\ast )}$ is empty. Note that this is equivalent to asserting that one
of the four inclusions $X\subset gX$, $X\subset gX^{\ast }$, $X^{\ast
}\subset gX$, $X^{\ast }\subset gX^{\ast }$ holds.

The following terminology will be useful.

\begin{definition}
A collection $E$ of subsets of $G$ which are closed under complementation is
called \textsl{nested} if for any pair $U$ and $V$ of sets in the
collection, one of the four sets $U^{(\ast )}\cap V^{(\ast )}$ is empty. If
each element $U$ of $E$ is a $H_{U}$--almost invariant subset of $G$ for some
subgroup $H_{U}$ of $G$, we will say that $E$ is \textsl{almost nested}%
\textit{\ }if for any pair $U$ and $V$ of sets in the collection, one of the
four sets $U^{(\ast )}\cap V^{(\ast )}$ is small.
\end{definition}

The above discussion shows that the translates of $X$ and $X^{\ast }$ under
the left action of $G$ are nested.

Note that $X$ is not uniquely determined by the splitting. In both cases, we
made choices of transversals, but it is easy to see that $X$ is independent
of the choice of transversal. However, in the case when $G=A\ast _{H}B$, we
chose $X$ to consist of elements for which $a_{1}$ is non-trivial whereas we
could equally well have reversed the roles of $A$ and $B$. This would simply
replace $X$ by $X^{\ast }-H$. Also either of these sets could be replaced by
its complement. We will use the term \textit{standard almost invariant set}
for the images in $H\backslash G$ of any one of $X$, $X\cup H$, $X^{\ast }$, 
$X^{\ast }-H$. In the case when $G=A\ast _{H}$, reversing the roles of the
two inclusion maps of $H$ into $A$ also replaces $X$ by $X^{\ast }-H$. Again
we have four standard almost invariant sets which are the images in $%
H\backslash G$ of any one of $X$, $X\cup H$, $X^{\ast }$, $X^{\ast }-H$.
There is a subtle point here. In the amalgamated free product case, we use
the obvious isomorphism between $A\ast _{H}B$ and $B\ast _{H}A$. In the HNN
case, let us write $A\ast _{H,i,j}$ to denote the group $%
<A,t:t^{-1}i(h)t=j(h)>$. Then the correct isomorphism to use between $A\ast
_{H,i,j}$ and $A\ast _{H,j,i}$ is not the identity on $A$. Instead it sends $%
t$ to $t^{-1}$ and $A$ to $t^{-1}At$. In all cases, we have four standard
almost invariant subsets of $H\backslash G.$

\begin{definition}
\label{defnofintersectionnumberofsplittings}If a group $G$ has splittings
over subgroups $H$ and $K$, and if $P$ and $Q$ are standard almost invariant
subsets of $H\backslash G$ and $K\backslash G$ respectively associated to
these splittings, then the \textsl{intersection number} of this pair of
splittings of $G$ is the intersection number of $P$ and $Q.$
\end{definition}

\begin{remark}
\label{XgisequivalenttoX}As any two of the four standard almost invariant
subsets of $H\backslash G$ associated to a splitting of $G$ over $H$ are
almost equal or almost complementary, Remark \ref
{almostequalsetshavesameintersectionnumber}\ tells us that this definition
does not depend on the choice of standard almost invariant subsets $P$ and $%
Q $.

If $X$ and $Y$ denote the pre-images in $G$ of $P\;$and $Q$ respectively,
and if we conjugate the first splitting by $a$ and the second by $b$, then $%
X $ is replaced by $aXa^{-1}$ and $Y\;$is replaced by $bYb^{-1}$. Now $Xg$
is $H$--almost equal to $X$ and $Yg$ is $K$--almost equal to $Y$, because of
the general fact that for any subset $W$ of $G$ and any element $g$ of $G$,
the set $Wg$ lies in a $l$--neighbourhood of $W$, where $l$ equals the length
of $g$. This follows from the equations $d(wg,w)=d(g,e)=l$. It follows that
the intersection number of a pair of splittings is unchanged if we replace
them by conjugate splittings.
\end{remark}

Now we can state two easy results about the case of zero intersection
number. Recall that if $X$ is one of the standard $H$--almost
invariant subsets of\textit{\ }$G$ determined by a splitting of $G$ over $H$%
, then the set of translates of $X\;$and $X^{\ast }$ is nested. It follows
at once that the self-intersection number of $H\backslash X$ is zero. Also
if two splittings of $G$ over subgroups $H$ and $K$ are compatible, and if $%
X $ and $Y$ denote corresponding standard $H$--almost and $K$--almost
invariant subsets of $G$, then the set of all translates of $X$, $X^{\ast }$%
, $Y$, $Y^{\ast }$ is also nested, so that the intersection number of the
two splittings is zero. The next section is devoted to proving converses to
each of these statements.

Before going further, we need to say a little more about splittings. Recall
from the introduction that a group $G$ is said to split over a subgroup $H$
if $G$ is isomorphic to $A\ast _{H}$ or to $A\ast _{H}B$, with $A\neq H\neq
B $. We will need a precise definition of a splitting. We will say that a 
\textit{splitting} of $G$ consists either of proper subgroups $A$ and $B$ of 
$G$ and a subgroup $H$ of $A\cap B$ such that the natural map $A\ast
_{H}B\rightarrow G$ is an isomorphism, or it consists of a subgroup $A$ of $%
G $ and subgroups $H_{0}$ and $H_{1}$ of $A$ such that there is an element $%
t $ of $G$ which conjugates $H_{0}$ to $H_{1}$ and the natural map $A\ast
_{H}\rightarrow G$ is an isomorphism.

Recall also that a collection of $n$ splittings of a group $G$ is \emph{%
compatible} if $G$ can be expressed as the fundamental group of a graph of
groups with $n$ edges, such that, for each $i$, collapsing all edges but the 
$i$-th yields the $i$-th splitting of $G.$ We note that if a splitting of a
group $G$ over a subgroup $H$ is compatible with a conjugate of itself by
some element $g$ of $G$, then $g$ must lie in $H$. This follows from a
simple analysis of the possibilities. For example, if the splitting $G=A\ast
_{H}B$ is compatible with its conjugate by some $g\in G$, then $G$ is the
fundamental group of a graph of groups with two edges, which must be a tree,
such that collapsing one edge yields the first splitting and collapsing the
other yields its conjugate by $g$. This means that each of the two extreme
vertex groups of the tree must be one of $A$, $A^{g}$, $B$ or $B^{g}$, and
the same holds for the subgroup of $G$ generated by the two vertex groups of
an edge. Now it is easy to see that $A\subset A^{g}$ and $B^{g}\subset B$,
or the same inclusions hold with the roles of $A$ and $B$ reversed. In
either case it follows that $g$ lies in $H$ as claimed. The case when $%
G=A\ast _{H}$ is slightly different, but the conclusion is the same. This
leads us to the following idea of equivalence of two splittings. We will say
that two amalgamated free product splittings of $G$ are equivalent, if they
are obtained from the same choice of subgroups $A$, $B$ and $H$ of $G$. This
means that the splittings $A\ast _{H}B$ and $B\ast _{H}A$ of $G$ are
equivalent. Similarly, a splitting $A\ast _{H}$ of $G$ is equivalent to the
splitting obtained by interchanging the two subgroups $H_{0}$ and $H_{1}$ of 
$A$. Also we will say that any splitting of a group $G$ over a subgroup $H$
is equivalent to any conjugate by some element of $H$. Then the equivalence
relation on all splittings of $G$ which this generates is the idea of
equivalence which we will need. Stated in this language, we see that if two
splittings are compatible and conjugate, then they must be equivalent.

Note that two splittings of a group $G$ are equivalent if and only if they
are over the same subgroup $H$, and they have exactly the same four standard
almost invariant sets.

Next we need to recall the connection between splittings of groups and
actions on trees. Bass--Serre theory, \cite{Serre} or \cite{Serre2}, tells us
that if a group $G$ splits over a subgroup $H$, then $G$ acts without
inversions on a tree $T$, so that the quotient is a graph with a single edge
and the vertex stabilisers are conjugate to $A$ or $B$ and the edge
stabilisers are conjugate to $H$. In his important paper \cite
{Dunwoody:Acc.1}, Dunwoody gave a method for constructing such a $G$--tree
starting from the subset $X$ of $G$ defined above. The crucial property of $%
X $ which is needed for the construction is the nestedness of the set of
translates of $X$ under the left action of $G$. We recall Dunwoody's result:

\begin{theorem}
\label{Dunwoodytreeconstruction}Let $E$ be a partially ordered set equipped
with an involution $e\rightarrow \overline{e}$, where $e\neq \overline{e}$,
such that the following conditions hold:

\begin{enumerate}
\item  If $e$, $f\in E$ and $e\leq f$, then $\overline{f}\leq \overline{e}$.

\item  If $e$, $f\in E$, there are only finitely many $g\in E$ such that $%
e\leq g\leq f$.

\item  If $e$, $f\in E$, at least one of the four relations $e\leq f$, $%
e\leq \overline{f}$, $\overline{e}\leq f$, $\overline{e}\leq \overline{f}$
holds.

\item  If $e$, $f\in E$, one cannot have $e\leq f$ and $e\leq \overline{f}.$
\end{enumerate}

Then there is an abstract tree $T$ with edge set equal to $E$ such that the
order relation which $E$ induces on the edge set of $T$ is equal to the
order relation in which $e\leq f$ if and only if there is an oriented path
in $T$ which begins with $e$ and ends with $f.$
\end{theorem}

One applies this result to the set $E=\{gX,gX^{\ast }:g\in G\}$ with the
partial order given by inclusion and the involution by complementation.
There is a natural action of $G$ on $E$ and hence on the tree $T$. In most
cases, $G$ acts on $T$ without inversions and we can recover the original
decomposition from this action as follows. Let $e$ denote the edge of $T$
determined by $X$. Then $X$ can be described as the set $\{g:g\in G,ge<e$ or 
$g\overline{e}<e\}$. If the action of $G$ on $T$ has inversions, then the
original splitting must have been an amalgamated free product decomposition $%
G=A\ast _{H}B$, with $H$ of index $2$ in $A$. In this case, subdividing the
edges of $T$ yields a tree $T_{1}$ on which $G$ acts without inversions. If $%
e_{1}$ denotes the edge of $T_{1}$ contained in $e$ and containing the
terminal vertex of $e$, then $X$ can be described as the set $\{g:g\in
G,ge_{1}<e_{1}$ or $g\overline{e_{1}}<e_{1}\}$.

Now we will prove the following result. This implies part 2) of Remark \ref
{almostequalsetshavesameintersectionnumber}. We give the proof here because
the proof in \cite{Scott:Intersectionnumbers} is not complete, and we will
need to apply the methods of proof later in this paper.

\begin{lemma}
\label{finitenumberofdoublecosets}Let $G$ be a finitely generated group with
finitely generated subgroups $H$ and $K$, a non-trivial $H$--almost invariant
subset $X$ and a non-trivial $K$--almost invariant subset $Y$. Then $\{g\in
G:gX$ and $Y$ are not nested\} consists of a finite number of double cosets $%
KgH.$
\end{lemma}

\begin{proof}
Let $\Gamma $ denote the Cayley graph of $G$ with respect to some finite
generating set for $G$. Let $P$ denote the almost invariant subset $%
H\backslash X$ of $H\backslash G$ and let $Q$ denote the almost invariant
subset $K\backslash Y$ of $K\backslash G$. Recall from the start of this
section, that if we identify $P$ with the $0$--cochain on $H\backslash \Gamma 
$ whose support is $P$, then $P$ is an almost invariant subset of $%
H\backslash G$ if and only if $\delta P$ is finite. Thus $\delta P$ is a
finite collection of edges in $H\backslash \Gamma $ and similarly $\delta Q$
is a finite collection of edges in $K\backslash \Gamma $. Now let $C$ denote
a finite connected subgraph of $H\backslash \Gamma $ such that $C$ contains $%
\delta P$ and the natural map $\pi _{1}(C)\rightarrow H$ is onto, and let $E$
denote a finite connected subgraph of $K\backslash \Gamma $ such that $E$
contains $\delta Q$ and the natural map $\pi _{1}(E)\rightarrow K$ is onto.
Thus the pre-image $D$ of $C$ in $\Gamma $ is connected and contains $\delta
X$, and the pre-image $F$ of $E$ in $\Gamma $ is connected and contains $%
\delta Y.$ Let $\Delta $ denote a finite subgraph of $D$ which projects onto 
$C$, and let $\Phi $ denote a finite subgraph of $F$ which projects onto $E$%
. If $gD$ meets $F$, there must be elements $h$ and $k$ in $H$ and $K$ such
that $gh\Delta $ meets $k\Phi $. Now $\{\gamma \in G:\gamma \Delta $ meets $%
\Phi \}$ is finite, as $G$ acts freely on $\Gamma $. It follows that $\{g\in
G:gD$ meets $F\}$ consists of a finite number of double cosets $KgH.$

The result would now be trivial if $X$ and $Y$ were each the vertex set of a
connected subgraph of $\Gamma $. As this need not be the case, we need to
make a careful argument as in the proof of Lemma 5.10 of \cite
{Scott-Wall:Topological}. Consider $g$ in $G$ such that $gD$ and $F$ are
disjoint. We will show that $gX$ and $Y$ are nested. As $D$ is connected,
the vertex set of $gD$ must lie entirely in $Y$ or entirely in $Y^{\ast }.$
Suppose that the vertex set of $gD$ lies in $Y$. For a set $S$ of vertices
of $\Gamma $, let $\overline{S}$ denote the maximal subgraph of $\Gamma $
with vertex set equal to $S$. Each component $W$ of $\overline{X}$ and $%
\overline{X^{\ast }}$ contains a vertex of $D$. Hence $gW$ contains a vertex
of $gD$ and so must meet $Y$. If $gW$ also meets $Y^{\ast }$, then it must
meet $F$. But as $F$ is connected and disjoint from $gD$, it lies in a
single component $gW$. It follows that there is exactly one component $gW$
of $\overline{gX}$ and $\overline{gX^{\ast }}$ which meets $Y^{\ast }$, so
that we must have $gX\subset Y$ or $gX^{\ast }\subset Y$. Similarly, if $gD$
lies in $Y^{\ast }$, we will find that $gX\subset Y^{\ast }$ or $gX^{\ast
}\subset Y^{\ast }$. It follows that in either case $gX$ and $Y$ are nested
as required.
\end{proof}

In Theorem 2.2 of \cite{Scott:TorusTheorem}, Scott used Dunwoody's theorem
to prove a general splitting result in the context of surfaces in $3$%
--manifolds. We will use the ideas in his proof a great deal. The following
theorem is the natural generalisation of his result to our more general
context and will be needed in the proofs of Theorems \ref{disjointsplittings}
and \ref{splittingsexist}. The first part of the theorem directly
corresponds to the result proved in \cite{Scott:TorusTheorem}, and the
second part is a simple generalisation which will be needed later.

\begin{theorem}$\phantom{111}$
\label{algebrafromTorusTheorem}
\begin{enumerate}
\item  Let $H$ be a finitely generated subgroup of a finitely generated
group $G$. Let $X$ be a non-trivial $H$--almost invariant set in $G$ such
that $E=\{gX,gX^{\ast }:g\in G\}$ is almost nested and if two of the four
sets $X^{(\ast )}\cap gX^{(\ast )}$ are small, then at least one of them is
empty. Then $G$ splits over the stabilizer $H^{\prime }$ of $X$ and $%
H^{\prime }$ contains $H$ as a subgroup of finite index. Further, one of the 
$H^{\prime }$--almost invariant sets $Y$ determined by the splitting is $H$%
--almost equal to $X.$

\item  Let $H_{1},\ldots ,H_{k}$ be finitely generated subgroups of a
finitely generated group $G$. Let $X_{i}$, $1\leq i\leq k$, be a non-trivial 
$H_{i}$--almost invariant set in $G$ such that $E=\{gX_{i},gX_{i}^{\ast
}:1\leq i\leq k,g\in G\}$ is almost nested. Suppose further that, for any
pair of elements $U$ and $V$ of $E$, if two of the four sets $U^{(\ast
)}\cap V^{(\ast )}$ are small, then at least one of them is empty. Then $G$
can be expressed as the fundamental group of a graph of groups whose $i$-th
edge corresponds to a conjugate of a splitting of $G$ over the stabilizer $%
H_{i}^{\prime }$ of $X_{i}$, and $H_{i}^{\prime }$ contains $H_{i}$ as a
subgroup of finite index. Further, for each $i$, one of the $H_{i}^{\prime }$%
--almost invariant sets determined by the $i$-th splitting is $H_{i}$--almost
equal to $X_{i}$.
\end{enumerate}
\end{theorem}

Most of the arguments needed to prove this theorem are contained in the
proof of Theorem 2.2 of \cite{Scott:TorusTheorem}, but in the context of $3$%
--manifolds. We will present the proof of the first part of this theorem, and
then briefly discuss the proof of the second part. The idea in the first
part is to define a partial order on $E=\{gX,gX^{\ast }:g\in G\}$, which
coincides with inclusion whenever possible. Let $U$ and $V$ denote elements
of $E$. If $U\cap V^{\ast }$ is small, we want to define $U\leq V$. There is
a difficulty, which is what to do if $U$ and $V$ are distinct but $U\cap
V^{\ast }$ and $V\cap U^{\ast }$ are both small. However, the assumption in
the statement of Theorem \ref{algebrafromTorusTheorem} is that if two of the
four sets $U^{(\ast )}\cap V^{(\ast )}$ are small, then one of them is
empty. Thus, as in \cite{Scott:TorusTheorem}, we define $U\leq V$ if and
only if $U\cap V^{\ast }$ is empty or the only small set of the four. Note
that if $U\subset V$ then $U\leq V$. We will show that this definition
yields a partial order on $E.$

As usual, we let $\Gamma $ denote the Cayley graph of $G$ with respect to
some finite generating set. The distance between two points of $G$ is the
usual one of minimal edge path length. Our first step is the analogue of
Lemma 2.3 of \cite{Scott:TorusTheorem}.

\begin{lemma}
\label{smallimpliesinboundednbhd}$U\cap V^{\ast }$ is small if and only if
it lies in a bounded neighbourhood of each of $U,U^{\ast }$, $V$, $V^{\ast
}. $
\end{lemma}

\begin{proof}
As $U$ and $V$ are translates of $X$ or $X^{\ast }$, it suffices to prove
that $gX\cap X^{\ast }$ is small if and only if it lies in a bounded
neighbourhood of each of $X$, $X^{\ast }$, $gX$, $gX^{\ast }$. If $gX\cap
X^{\ast }$ is small, it projects to a finite subset of $H\backslash G$ which
therefore lies within a bounded neighbourhood of the image of $\delta X$. By
lifting paths, we see that each point of $gX\cap X^{\ast }$ lies in a
bounded neighbourhood of $\delta X$, and hence lies in a bounded
neighbourhood of $X$ and $X^{\ast }$. By reversing the roles of $gX$ and $%
X^{\ast }$, we also see that $gX\cap X^{\ast }$ lies in a bounded
neighbourhood of each of $gX$ and $gX^{\ast }.$

For the converse, suppose that $gX\cap X^{\ast }$ lies in a bounded
neighbourhood of each of $X$ and $X^{\ast }$. Then it must lie in a bounded
neighbourhood of $\delta X$, so that its image in $H\backslash G$ must lie
in a bounded neighbourhood of the image of $\delta X$. As this image is
finite, it follows that $gX\cap X^{\ast }$ must be small, as required.
\end{proof}

Now we can prove that our definition of $\leq $ yields a partial order on $%
E. $ Our proof is essentially the same as in Lemma 2.4 of \cite
{Scott:TorusTheorem}.

\begin{lemma}
\label{partialorder}If a relation $\leq $ is defined on $E$ by the condition
that $U\leq V$ if and only if $U\cap V^{\ast }$ is empty or the only small
set of the four sets $U^{(\ast )}\cap V^{(\ast )}$, then $\leq $ is a
partial order.
\end{lemma}

\begin{proof}
We need to show that $\leq $ is transitive and that if $U\leq V$ and $V\leq
U $ then $U=V.$

Suppose first that $U\leq V$ and $V\leq U$. The first inequality implies
that $U\cap V^{\ast }$ is small and the second implies that $V\cap U^{\ast }$
is small, so that two of the four sets $U^{(\ast )}\cap V^{(\ast )}$ are
small. The assumption of Theorem \ref{algebrafromTorusTheorem} implies that
one of these two sets must be empty. As $U\leq V$, our definition of $\leq $
implies that $U\cap V^{\ast }$ is empty. Similarly, the fact that $V\leq U$
tells us that $V\cap U^{\ast }$ is empty. This implies that $U=V$ as
required.

To prove transitivity, let $U$, $V$ and $W$ be elements of $E$ such that $%
U\leq V\leq W$. We must show that $U\leq W.$

Our first step is to show that $U\cap W^{\ast }$ is small. As $U\cap V^{\ast
}$ and $V\cap W^{\ast }$ are small, we let $d_{1}$ be an upper bound for the
distance of points of $U\cap V^{\ast }$ from $V$ and let $d_{2}$ be an upper
bound for the distance of points of $V\cap W^{\ast }$ from $W$. Let $x$ be a
point of $U\cap W^{\ast }$. If $x$ lies in $V$, then it lies in $V\cap
W^{\ast }$ and so has distance at most $d_{2}$ from $W$. Otherwise, it must
lie in $U\cap V^{\ast }$ and so have distance at most $d_{1}$ from some
point $x^{\prime }$ of $V$. If $x^{\prime }$ lies in $W$, then $x$ has
distance at most $d_{1}$ from $W$. Otherwise, $x^{\prime }$ lies in $V\cap
W^{\ast }$ and so has distance at most $d_{2}$ from $W$. In this case, $x$
has distance at most $d_{1}+d_{2}$ from $W$. It follows that in all cases, $%
x $ has distance at most $d_{1}+d_{2}$ from $W$, so that $U\cap W^{\ast }$
lies in a bounded neighbourhood of $W$ as required. As $U\cap W^{\ast }$ is
contained in $W^{\ast }$, it follows that it lies in bounded neighbourhoods
of $W$ and $W^{\ast }$, so that $U\cap W^{\ast }$ is small as required.

The definition of $\leq $ now shows that $U\leq W$, except possibly when two
of the four sets $U^{(\ast )}\cap W^{(\ast )}$ are small. The only
possibility is that $U^{\ast }\cap W$ and $U\cap W^{\ast }$ are both small.
As one must be empty, either $U\subset W$ or $W\subset U$. We conclude that
if $U\leq V\leq W$, then either $U\leq W$ or $W\subset U$. Now we consider
two cases.

First suppose that $U\subset V\leq W$, so that either $U\leq W$ or $W\subset
U$. If $W\subset U$, then $W\subset V$, so that $W\leq V$. As $V\leq W$ and $%
W\leq V$, it follows from the first paragraph of the proof of this lemma
that $V=W$. Hence, in either case, $U\leq W.$

Now consider the general situation when $U\leq V\leq W$. Again either $U\leq
W$ or $W\subset U$. If $W\subset U$, then we have $W\subset U\leq V$. Now
the preceding paragraph implies that $W\leq V$. Hence we again have $V\leq W$
and $W\leq V$ so that $V=W$. Hence $U\leq W$ still holds. This completes the
proof of the lemma.
\end{proof}

Next we need to verify that the set $E$ with the partial order which we have
defined satisfies all the hypotheses of Dunwoody's Theorem \ref
{Dunwoodytreeconstruction}.

\begin{lemma}
\label{posatisfiesDunwoody}$E$ together with $\leq $ satisfies the following
conditions.

\begin{enumerate}
\item  If $U$, $V\in E$ and $U\leq V$, then $V^{\ast }\leq U^{\ast }$.

\item  If $U$, $V\in E$, there are only finitely many $Z\in E$ such that $%
U\leq Z\leq V$.

\item  If $U$, $V\in E$, at least one of the four relations $U\leq V$, $%
U\leq V^{\ast }$, $U^{\ast }\leq V$, $U^{\ast }\leq V^{\ast }$ holds.

\item  If $U$, $V\in E$, one cannot have $U\leq V$ and $U\leq V^{\ast }.$
\end{enumerate}
\end{lemma}

\begin{proof}
Conditions (1) and (3) are obvious from the definition of $\leq $ and the
hypotheses of Theorem \ref{algebrafromTorusTheorem}.

To prove (4), we observe that if $U\leq V$ and $U\leq V^{\ast }$, then $%
U\cap V^{\ast }$ and $U\cap V$ must both be small. This implies that $U$
itself is small, so that $X$ or $X^{\ast }$ must be small. But this
contradicts the hypothesis that $X$ is a non-trivial $H$--almost invariant
subset of $G.$

Finally we prove condition (2). Let $Z=gX$ be an element of $E$ such that $%
Z\leq X$. Recall that, as $Z\cap X^{\ast }$ projects to a finite subset of $%
H\backslash G$, we know that $Z\cap X^{\ast }$ lies in a $d$--neighbourhood
of $X$, for some $d>0$. If $Z\leq X$ but $Z$ is not contained in $X$, then $%
Z $ and $X$ are not nested. Now Lemma \ref{finitenumberofdoublecosets} tells
us that if $Z$ is such a set, then $g$ belongs to one of only finitely many
double cosets $HkH$. It follows that if we consider all elements $Z$ of $E$
such that $Z\leq X$, we will find either $Z\subset X$, or $Z\cap X^{\ast }$
lies in a $d$--neighbourhood of $X$, for finitely many different values of $%
d. $ Hence there is $d_{1}>0$ such that if $Z\leq X$ then $Z$ lies in the $%
d_{1} $--neighbourhood of $X$. Similarly, there is $d_{2}>0$ such that if $%
Z\leq X^{\ast }$, then $Z$ lies in the $d_{2}$--neighbourhood of $X^{\ast }.$
Let $d $ denote the larger of $d_{1}$ and $d_{2}$. Then for any elements $U$
and $V$ of $E$ with $U\leq V$, the set $U\cap V^{\ast }$ lies in the $d$%
--neighbourhood of each of $U$, $U^{\ast }$, $V$ and $V^{\ast }.$

Now suppose we are given $U\leq V$ and wish to prove condition (2). Choose a
point $u$ in $U$ whose distance from $U^{\ast }$ is greater than $d$, choose
a point $v$ in $V^{\ast }$ whose distance from $V$ is greater than $d$ and
choose a path $L$ in $\Gamma $ joining $u$ to $v$. If $U\leq Z\leq V$, then $%
u$ must lie in $Z$ and $v$ must lie in $Z^{\ast }$ so that $L$ must meet $%
\delta Z$. As $L$ is compact, the proof of Lemma \ref
{finitenumberofdoublecosets} shows that the number of such $Z$ is finite.
This completes the proof of part 2) of the lemma.
\end{proof}

\medskip

We are now in a position to prove Theorem \ref{algebrafromTorusTheorem}.

\begin{proof}
To prove the first part, we let $E$ denote the set of all translates of $X$
and $X^{\ast }$ by elements of $G$, let $U\rightarrow U^{\ast }$ be the
involution on $E$ and let the relation $\leq $ be defined on $E$ by the
condition that $U\leq V$ if $U\cap V^{\ast }$ is empty or the only small set
of the four sets $U^{(\ast )}\cap V^{(\ast )}$. Lemmas \ref{partialorder}
and \ref{posatisfiesDunwoody} show that $\leq $ is a partial order on $E$
and satisfies all of Dunwoody's conditions (1)--(4). Hence we can construct a
tree $T$ from $E$. As $G$ acts on $E$, we have a natural action of $G$ on $%
T. $ Clearly, $G$ acts transitively on the edges of $T$. If $G$ acts without
inversions, then $G\backslash T$ has a single edge and gives $G$ the
structure of an amalgamated free product or HNN decomposition. The
stabiliser of the edge of $T$ which corresponds to $X$ is the stabiliser $%
H^{\prime }$ of $X$, so we obtain a splitting of $G$ over $H^{\prime }$
unless $G$ fixes a vertex of $T$. Note that as $H\backslash \delta X$ is
finite, and $H^{\prime }$ preserves $\delta X$, it follows that $H^{\prime }$
contains $H$ with finite index as claimed in the theorem. If $G$ acts on $T$
with inversions, we simply subdivide each edge to obtain a new tree $%
T^{\prime }$ on which $G$ acts without inversions. In this case, the
quotient $G\backslash T^{\prime }$ again has one edge, but it has distinct
vertices. The edge group is $H^{\prime }$ and one of the vertex groups
contains $H^{\prime }$ with index two. As $H$ has infinite index in $G$, it
follows that in this case also we obtain a splitting of $G$ unless $G$ fixes
a vertex of $T.$

Suppose that $G$ fixes a vertex $v$ of $T$. As $G$ acts transitively on the
edges of $T$, every edge of $T$ must have one vertex at $v$, so that all
edges of $T$ are adjacent to each other. We will show that this cannot
occur. The key hypothesis here is that $X$ is non-trivial.

Let $W$ denote $\{g:gX\leq X$ or $gX^{\ast }\leq X\}$, and note that
condition 3) of Lemma \ref{posatisfiesDunwoody} shows that $W^{\ast
}=\{g:gX\leq X^{\ast }$ or $gX^{\ast }\leq X^{\ast }\}$. Recall that there
is $d_{1}>0$ such that if $Z\leq X$ then $Z$ lies in the $d_{1}$%
--neighbourhood of $X$. If $d$ denotes $d_{1}+1$, and $g\in W$, it follows
that $g\delta X$ lies in the $d$--neighbourhood of $X$. Let $c$ denote the
distance of the identity of $G$ from $\delta X$. Then $g$ must lie within
the $(c+d)$--neighbourhood of $X$, for all $g\in W$, so that $W$ itself lies
in the $(c+d)$--neighbourhood of $X$. Similarly, $W^{\ast }$ lies in the $%
(c+d)$--neighbourhood of $X^{\ast }$. Now both $X$ and $X^{\ast }$ project to
infinite subsets of $H\backslash G$, so $G$ cannot equal $W$ or $W^{\ast }.$
It follows that there are elements $U\;$and $V$ of $E$ such that $U<X<V$, so
that $U\;$and $V$ represent non-adjacent edges of $T$. This completes the
proof that $G$ cannot fix a vertex of $T.$

To prove the last statement of the first part of Theorem \ref
{algebrafromTorusTheorem}, we will simplify notation by supposing that the
stabiliser $H^{\prime }$ of $X$ is equal to $H$. One of the standard $H$%
--almost invariant sets associated to the splitting we have obtained from the
action of $G$ on the tree $T$ is the set $W$ in the preceding paragraph. We
will show that $W$ is $H$--almost equal to $X$. The preceding paragraph shows
that $W$ lies in the $(c+d)$--neighbourhood of $X$, and that $W^{\ast }$ lies
in the $(c+d)$--neighbourhood of $X^{\ast }$. It follows that $W$ is $H$%
--almost contained in $X$ and $W^{\ast }$ is $H$--almost contained in $X^{\ast
}$, so that $W$ and $X$ are $H$--almost equal as claimed. This completes the
proof of the first part of Theorem \ref{algebrafromTorusTheorem}.

For the second part, we will simply comment on the modifications needed to
the preceding proof. The statement of Lemma \ref{smallimpliesinboundednbhd}
remains true though the proof needs a little modification. The statement and
proof of Lemma \ref{partialorder} apply unchanged. The statement of Lemma 
\ref{posatisfiesDunwoody} remains true, though the proof needs some minor
modifications. Finally the proof of the first part of Theorem \ref
{algebrafromTorusTheorem} applies with minor modifications to show that $G$
acts on a tree $T$ with quotient consisting of $k$ edges in the required
way. This completes the proof of Theorem \ref{algebrafromTorusTheorem}.
\end{proof}

\section{Zero intersection numbers}

In this section, we prove our two main results about the case of zero
intersection number. First we will need the following little result.

\begin{lemma}
\label{bignormaliserimpliesnormal}Let $G$ be a finitely generated group
which splits over a subgroup $H$. If the normaliser $N$ of $H$ in $G$ has
finite index in $G$, then $H$ is normal in $G.$
\end{lemma}

\begin{proof}
The given splitting of $G$ over $H$ corresponds to an action of $G$ on a
tree $T$ such that $G\backslash T\;$has a single edge, and some edge of $T\;$%
has stabiliser $H$. Let $T^{\prime }$ denote the fixed set of $H$, ie, the
set of all points fixed by $H$. Then $T^{\prime }$ is a (non-empty) subtree
of $T$. As $N$ normalises $H$, it must preserve $T^{\prime }$, ie $%
NT^{\prime }=T^{\prime }$. Suppose that $N\neq G$. As $N$ has finite index
in $G$, we let $e,g_{1},\ldots ,g_{n}$ denote a set of coset representatives
for $N$ in $G$, where $n\geq 1$. As $G$ acts transitively on $T$, we have $%
T=T^{\prime }\cup g_{1}T^{\prime }\cup \ldots \cup g_{n}T^{\prime }$. Edges
of $T^{\prime }$ all have stabiliser $H$, and so edges of $g_{i}T^{\prime }$
all have stabiliser $g_{i}Hg_{i}^{-1}$. As $g_{i}$ does not lie in $N$,
these stabilisers are distinct so the intersection $T^{\prime }\cap
g_{i}T^{\prime }$ contains no edges. The intersection of two subtrees of a
tree must be empty or a tree, so it follows that $T^{\prime }\cap
g_{i}T^{\prime }$ is empty or a single vertex $v_{i}$, for each $i$. Now $N$
preserves $T^{\prime }$ and permutes the translates $g_{i}T^{\prime }$, so $%
N $ preserves the collection of all the $v_{i}$'s. As this collection is
finite, $N$ has a subgroup $N_{1}$ of finite index such that $N_{1}$ fixes a
vertex $v$ of $T^{\prime }$. As $N_{1}$ has finite index in $G$, it follows
that $G$ itself fixes some vertex of $T$, which contradicts our assumption
that our action of $G$ on $T$ corresponds to a splitting of $G$. This
contradiction shows that $N$ must equal $G$, so that $H$ is normal in $G\;$%
as claimed.
\end{proof}

Recall that if $X$ is a $H$--almost invariant subset of $G$ associated to a
splitting of $G$, then the set of translates of $X$ and $X^{\ast }$ is
nested. Equivalently, for every $g\in G$, one of the four sets $X^{(\ast
)}\cap gX^{(\ast )}$ is empty. We need to consider carefully how it is
possible for two of the four sets to be small, and a similar question arises
when one considers two splittings of $G.$

\begin{lemma}
\label{twosmallsets}Let $G$ be a finitely generated group with two
splittings over finitely generated subgroups $H$ and $K$ with associated $H$%
--almost invariant subset $X$ of $G$ and associated $K$--almost invariant
subset $Y$ of $G.$

\begin{enumerate}
\item  If two of the four sets $X^{(\ast )}\cap Y^{(\ast )}$ are small, then 
$H=K.$

\item  If two of the four sets $X^{(\ast )}\cap gX^{(\ast )}$ are small,
then $g$ normalises $H.$
\end{enumerate}
\end{lemma}

\begin{proof}
Our first step will be to show that $H$ and $K$ must be commensurable.
Without loss of generality, we can suppose that $X\cap Y$ is small. The
other small set can only be $X^{\ast }\cap Y^{\ast }$, as otherwise $X$ or $%
Y $ would be small which is impossible. It follows that for each edge of $%
\delta Y$, either it is also an edge of $\delta X$ or it has (at least) one
end in one of the two small sets. As the images in $H\backslash \Gamma $ of $%
\delta X$ and of each small set is finite, and as the graph $\Gamma $ is
locally finite, it follows that the image of $\delta Y$ in $H\backslash
\Gamma $ must be finite. This implies that $H\cap K$ has finite index in the
stabiliser $K$ of $\delta Y$. By reversing the roles of $H$ and $K$, it
follows that $H\cap K$ has finite index in $H$, so that $H$ and $K$ must be
commensurable, as claimed.

Now let $L$ denote $H\cap K$, so that $L$ stabilises both $X$ and $Y$, and
consider the images $P$ and $Q$ of $X$ and $Y$ in $L\backslash \Gamma $. As $%
L$ has finite index in $H$ and $K$, it follows that $\delta P$ and $\delta Q$
are each finite, so that $P$ and $Q$ are almost invariant subsets of $%
L\backslash G$. Further, two of the four sets $X^{(\ast )}\cap Y^{(\ast )}$
have finite image in $L\backslash \Gamma $, so we can assume that $P$ and $Q$
are almost equal, by replacing one of $X$ or $Y$ by its complement in $G$,
if needed. Let $L^{\prime }$ denote the intersection of the conjugates of $L$
in $H$, so that $L^{\prime }$ is normal in $H$, though it need not be normal
in $K$. We do not have $L^{\prime }=H\cap K$, but because $L$ has finite
index in $H$, we know that $L^{\prime }$ has finite index in $H$ and hence
also in $K$, which is all we need. Let $P^{\prime }$ and $Q^{\prime }$
denote the images of $X$ and $Y$ respectively in $L^{\prime }\backslash
\Gamma $, and consider the action of an element $h$ of $H$ on $L^{\prime
}\backslash \Gamma .$ Trivially $hP^{\prime }=P^{\prime }$. As $P^{\prime }$
and $Q^{\prime }$ are almost equal, $hQ^{\prime }$ must be almost equal to $%
Q^{\prime }$. Now we use the key fact that $Y$ is associated to a splitting
of $G$ so that its translates by $G$ are nested. Thus for any element $g$ of 
$G$, one of the following four inclusions holds: $gY\subset Y$, $gY\subset
Y^{\ast }$, $gY^{\ast }\subset Y$, $gY^{\ast }\subset Y^{\ast }$. As $%
hQ^{\prime }$ is almost equal to $Q^{\prime }$, we must have $hY\subset Y$
or $hY^{\ast }\subset Y^{\ast }.$ But $h$ has a power which lies in $L$ and
hence stabilises $Y$. It follows that $hY=Y$, so that $h$ lies in $K$. Thus $%
H$ is a subgroup of $K.$ Similarly, $K$ must be a subgroup of $H$, so that $%
H=K$. This completes the proof of part 1 of the lemma. Note that it follows
that $L=H=K$, that $H\backslash X=P$ and $K\backslash Y=Q$ and that $P$ and $%
Q$ are almost equal or almost complementary.

In order to prove part 2 of the lemma, we apply the preceding work to the
case when the second splitting is obtained from the first by conjugating by
some element $g$ of $G$. Thus $K=gHg^{-1}$ and $Y=gXg^{-1}$ which is $K$%
--almost equal to $gX$ by Remark \ref{XgisequivalenttoX}. Hence if two of the
four sets $X^{(\ast )}\cap gX^{(\ast )}$ are small, then so are two of the
four sets $X^{(\ast )}\cap Y^{(\ast )}$ small. Now the above shows that $%
H=K=gHg^{-1}$, so that $g$ normalises $H$. This completes the proof of the
lemma.
\end{proof}

\begin{lemma}
\label{equalsplittings}Let $G$ be a finitely generated group with two
splittings over finitely generated subgroups $H$ and $K$ with associated $H$%
--almost invariant subset $X$ of $G$ and associated $K$--almost invariant
subset $Y$ of $G$. If two of the four sets $X^{(\ast )}\cap Y^{(\ast )}$ are
small, then the two splittings of $G$ are conjugate. Further one of the
following holds:

\begin{enumerate}
\item  the two splittings are equivalent, or

\item  the two splittings are of the form $G=L\ast _{H}C$, where $H$ has
index $2$ in $L$, and the splittings are conjugate by an element of $L$, or

\item  $H$ is normal in $G$ and $H\backslash G$ is isomorphic to $\mathbb{Z}$
or to $\mathbb{Z}_{2}\ast \mathbb{Z}_{2}.$
\end{enumerate}
\end{lemma}

\begin{proof}
The preceding lemma showed that the hypotheses imply that $H$ equals $K$ and
also that the images $P$ and $Q$ of $X$ and $Y$ in $H\backslash G$ are
almost equal or almost complementary. By replacing one of $X$ or $Y$ by its
complement if needed, we can arrange that $P$ and $Q$ are almost equal. We
will show that in most cases, the two given splittings over $H$ and $K$ must
be equivalent, and that the exceptional cases can be analysed separately to
show that the splittings are conjugate.

Recall that by applying Theorem \ref{Dunwoodytreeconstruction}, we can use
information about $X$ and its translates to construct a $G$--tree $T_{X}$ and
hence the original splitting of $G$ over $H$. Similarly, we can use
information about $Y$ and its translates to construct a $G$--tree $T_{Y}$ and
hence the original splitting of $G$ over $K$. We will compare these two
constructions in order to prove our result.

As $P$ and $Q$ are almost equal subsets of $H\backslash G$, it follows that
there is $\delta \geq 0$ such that, in the Cayley graph $\Gamma $ of $G$, we
have $X$ lies in a $\delta $--neighbourhood of $Y$ and $Y$ lies in a $\delta $%
--neighbourhood of $X$. Now let $U_{X}$ denote one of $X$ or $X^{\ast }$, let 
$V_{X}$ denote one of $gX$ or $gX^{\ast }$ and let $U_{Y}$ and $V_{Y}$
denote the corresponding sets obtained by replacing $X$ with $Y$. Recall
that $U_{X}\cap V_{X}$ is small if and only if its image in $H\backslash G$
is finite. Clearly this occurs if and only if $V_{X}$ lies in a $\delta $%
--neighbourhood of $U_{X}^{\ast }$, for some $\delta \geq 0$. It follows that 
$U_{X}\cap V_{X}$ is small if and only if $U_{Y}\cap V_{Y}$ is small.

As $X$ and $Y$ are associated to splittings, we know that for each $g\in G$,
at least one of the four sets $X^{(\ast )}\cap gX^{(\ast )}$ is empty and at
least one of the four sets $Y^{(\ast )}\cap gY^{(\ast )}$ is empty. Further
the information about which of the four sets is empty completely determines
the trees $T_{X}$ and $T_{Y}$. Thus we would like to show that when we
compare the four sets $X^{(\ast )}\cap gX^{(\ast )}$ with the four sets $%
Y^{(\ast )}\cap gY^{(\ast )}$, then corresponding sets are empty. Note that
when $g$ lies in $H$, we have $gX=X$, so that two of the four sets $X^{(\ast
)}\cap gX^{(\ast )}$ are empty.

First we consider the case when, for each $g\in G-H$, only one of the sets $%
X^{(\ast )}\cap gX^{(\ast )}$ is small and hence empty. Then only the
corresponding one of the four sets $Y^{(\ast )}\cap gY^{(\ast )}$ is small
and hence empty. Now the correspondence $gX\rightarrow gY$ gives a $G$%
--isomorphism of $T_{X}$ with $T_{Y}$ and thus the splittings are equivalent.

Next we consider the case when two of the sets $X^{(\ast )}\cap gX^{(\ast )}$
are small, for some $g\in G-H$. Part 2 of Lemma \ref{twosmallsets} implies
that $g$ normalises $H$. Further if $R=H\backslash gX$, then $P\;$is almost
equal to $R$ or $R^{\ast }$. Let $N(H)$ denote the normaliser of $H$ in $G$,
so that $N(H)$ acts on the left on the graph $H\backslash \Gamma $ and we
have $R=gP$. Let $L$ denote the subgroup of $N(H)$ consisting of elements $k$
such that $kP$ is almost equal to $P$ or $P^{\ast }$. Now we apply Theorem
5.8 from \cite{Scott-Wall:Topological} to the action of $H\backslash L$ on
the left on the graph $H\backslash \Gamma .$ This result tells us that if $%
H\backslash L$ is infinite, then it has an infinite cyclic subgroup of
finite index. Further the proof of this result in \cite
{Scott-Wall:Topological} shows that the quotient of $H\backslash \Gamma $ by 
$H\backslash L$ must be finite. This implies that $H\backslash \Gamma $ has
two ends and that $L$ has finite index in $G$. To summarise, either $%
H\backslash L$ is finite, or it has two ends and $L$ has finite index in $G$%
. Let $k$ be an element of $L$ whose image in $H\backslash L$ has finite
order such that $kP\overset{a}{=}P$. As $X$ is associated to a splitting of $%
G$, we must have $kX\subset X$ or $X\subset kX$. As $k$ has finite order in $%
H\backslash L$, we have $k^{n}X=X$, for some positive integer $n$, which
implies that $kX=X$ so that $k$ itself lies in $H$. It follows that the
group $H\backslash L$ must be trivial, $\mathbb{Z}_{2}$, $\mathbb{Z}$ or $%
\mathbb{Z}_{2}\ast \mathbb{Z}_{2}$. In the first case, the two trees $T_{X}$
and $T_{Y}$ will be $G$--isomorphic, showing that the given splittings are
equivalent. In the other three cases, $L-H$ is non-empty and we know that,
for any $g\in L-H$, two of the four sets $X^{(\ast )}\cap gX^{(\ast )}$ are
small. Thus in these cases, it seems possible that $T_{X}$ and $T_{Y}$ will
not be $G$--isomorphic, so we need some special arguments.

We start with the case when $H\backslash L$ is $\mathbb{Z}_{2}$. In\ this
case, the given splitting must be an amalgamated free product of the form $%
L\ast _{H}C$, for some group $C$. If $k$ denotes an element of $L-H$, then $%
kP\overset{a}{=}P^{\ast }$. Thus $G$ acts on $T_{X}$ and $T_{Y}$ with
inversions. Recall that either the two partial orders on the translates of $%
X $ and $Y$ are the same under the bijection $gX\rightarrow gY$, or they
differ only in that $kX^{\ast }\subset X$ but $Y\subset kY^{\ast }$, for all 
$k\in L-H$. If they differ, we replace the second splitting by its conjugate
by some element $k\in L-H$, so that $Y$ is replaced by $Y^{\prime }=kY$ and
we replace $X$ by $X^{\prime }=X^{\ast }.$ As $Y^{\prime }$ is $H$--almost
equal to $X^{\prime }$, the partial orders on the translates of $X^{\prime }$
and $Y^{\prime }$ respectively are the same under the bijection $gX^{\prime
}\rightarrow gY^{\prime }$ except possibly when one compares $X^{\prime }$, $%
kX^{\prime }$ and $Y^{\prime }$, $kY^{\prime }$, where $k\in L-H.$ In this
case, the inclusion $kX^{\ast }\subset X$ tells us that $kX^{\prime }\subset
(X^{\prime })^{\ast }$, and the inclusion $Y\subset kY^{\ast }$ tells us
that $kY^{\prime }=k^{2}Y=Y\subset kY^{\ast }=\left( Y^{\prime }\right)
^{\ast }$. We conclude that the partial orders on the translates of $%
X^{\prime }$ and $Y^{\prime }$ respectively are exactly the same, so that $%
T_{X}$ and $T_{Y}$ are $G$--isomorphic, and the two given splittings are
conjugate by an element of $L$.

Now we turn to the two cases where $H\backslash L$ is infinite, so that $L$
has finite index in $G$ and $H\backslash \Gamma $ has two ends. As $L$
normalises $H$, Lemma \ref{bignormaliserimpliesnormal} shows that $H$ is
normal in $G$. As $H\backslash \Gamma $ has two ends, it follows that $L=G$,
so that $H\backslash G$ is $\mathbb{Z}$ or $\mathbb{Z}_{2}\ast \mathbb{Z}%
_{2} $. It is easy to check that there is only one splitting of $\mathbb{Z}$
over the trivial group and that all splittings of $\mathbb{Z}_{2}\ast 
\mathbb{Z}_{2}$ over the trivial group are conjugate. It follows that, in
either case, all splittings of $G$ over $H$ are conjugate. This completes
the proof of Lemma \ref{equalsplittings}.
\end{proof}

\begin{lemma}
\label{emptyorconjugate}Let $G$ be a finitely generated group with two
splittings over finitely generated subgroups $H$ and $K$ with associated $H$%
--almost invariant subset $X$ of $G$ and associated $K$--almost invariant
subset $Y$ of $G$. Let $E=\{gX,gX^{\ast },gY,gY^{\ast }:g\in G\}$, and let $%
U $ and $V$ denote two elements of $E$ such that two of the four sets $%
U^{(\ast )}\cap V^{(\ast )}$ are small. Then either one of the two sets is
empty, or the two given splittings of $G$ are conjugate.
\end{lemma}

\begin{proof}
Recall that $X$ is associated to a splitting of $G$ over $H$. It follows
that $gX$ is associated to the conjugate of this splitting by $g$. Thus $U$
and $V$ are associated to splittings of $G$ which are each conjugate to one
of the two given splittings. If $U$ and $V$ are each translates of $X$ or $%
X^{\ast }$, the nestedness of the translates of $X$ shows that one of the
two small sets must be empty as claimed. Similarly if both are translates of 
$Y$ or $Y^{\ast }$, then one of the two small sets must be empty. If $U$ is
a translate of $X$ or $X^{\ast }$ and $V$ is a translate of $Y$ or $Y^{\ast
} $, we apply Lemma \ref{equalsplittings} to show that the splittings to
which $U$ and $V$ are associated are conjugate. It follows that the two
original splittings were conjugate as required.
\end{proof}

Now we come to the proof of our first main result.

\begin{theorem}
\label{disjointsplittings}{Let }$G${\ be a finitely generated
group with }$n${\ splittings over finitely generated subgroups. This
collection of splittings is compatible up to conjugacy if and only if each
pair of splittings has intersection number zero. Further, in this situation,
the graph of groups structure on }$G$ obtained from these splittings has a
unique underlying graph, and the edge and vertex groups are unique up to
conjugacy.
\end{theorem}

\begin{proof}
Let the $n$ splittings $s_{i}$ of $G$ be over subgroups $H_{1},...,H_{n}$
with associated $H_{i}$--almost invariant subsets $X_{i}$ of $G$, and let $%
E=\{gX_{i},gX_{i}^{\ast }:g\in G,1\leq i\leq n\}$. We will start by
supposing that no two of the $s_{i}$'s are conjugate. We will handle the
general case at the end of this proof.

We will apply the second part of Theorem \ref{algebrafromTorusTheorem} to $%
E. $ Recall that our assumption that the $s_{i}$'s have intersection number
zero implies that no translate of $X_{i}$ can cross any translate of $X_{j}$%
, for $1\leq i\neq j\leq n$. As each $X_{i}$ is associated to a splitting,
it is also true that no translate of $X_{i}$ can cross any translate of $%
X_{i}$. This means that the set $E$ is almost nested. In order to apply
Theorem \ref{algebrafromTorusTheorem}, we will also need to show that for
any pair of elements $U$ and $V$ of $E$, if two of the four sets $U^{(\ast
)}\cap V^{(\ast )}$ are small then one is empty. Now Lemma \ref
{emptyorconjugate} shows that if two of these four sets are small, then
either one is empty or there are distinct $i$ and $j$ such that $s_{i}$ and $%
s_{j}$ are conjugate. As we are assuming that no two of these splittings are
conjugate, it follows that if two of the four sets $U^{(\ast )}\cap V^{(\ast
)}$ are small then one is empty, as required.

Theorem \ref{algebrafromTorusTheorem} now implies that $G$ can be expressed
as the fundamental group of a graph $\Gamma $ of groups whose $i$-th edge
corresponds to a conjugate of a splitting of $G$ over the stabilizer $%
H_{i}^{\prime }$ of $X_{i}$. As $X_{i}$ is associated to a splitting of $G$
over $H_{i}$, its stabiliser $H_{i}^{\prime }$ must equal $H_{i}$. Further,
it is clear from the construction that collapsing all but the $i$-th edge of 
$\Gamma $ yields a conjugate of $s_{i}$, as the corresponding $G$--tree has
edges which correspond precisely to the translates of $X_{i}.$

Now suppose that we have a graph of groups structure $\Gamma ^{\prime }$ for 
$G$ such that, for each $i$, $1\leq i\leq n$, collapsing all edges but the $%
i $-th yields a conjugate of the splitting $s_{i}$ of $G$. This determines
an action of $G$ on a tree $T^{\prime }$ without inversions. We want to show
that $T$ and $T^{\prime }$ are $G$--isomorphic. For this implies that $\Gamma 
$ and $\Gamma ^{\prime }$ have the same underlying graph, and that
corresponding edge and vertex groups are conjugate, as required. Let $e$
denote an edge of $T^{\prime }$, and let $Y(e)$ denote $\{g\in G:ge<e$ or $g%
\overline{e}<e\}$. There are edges $e_{i}$ of $T^{\prime }$, $1\leq i\leq n$%
, such that the set $E^{\prime }$ of all translates of $Y(e_{i})$ and $%
Y(e_{i})^{\ast }$ is nested and Dunwoody's construction applied to $%
E^{\prime }$ yields the $G$--tree $T^{\prime }$ again. We will denote $%
Y(e_{i})$ by $Y_{i}$. The hypotheses imply that there is $k\in G$ such that
the stabiliser $K_{i}$ of $e_{i}$ equals $k^{-1}H_{i}k$, and that $Y_{i}$ is 
$K_{i}$--almost equal to $k^{-1}X_{i}k$, where $X_{i}$ is one of the standard 
$H_{i}$--almost invariant subsets of $G$ associated to the splitting $s_{i}$.
Let $Z_{i}$ denote $kY_{i}$ so that $Z_{i}$ is $H_{i}$--almost equal to $%
X_{i}k$. Now Remark \ref{XgisequivalenttoX} shows that $X_{i}k$ is $H_{i}$%
--almost equal to $X_{i}$, so that $Z_{i}$ is $H_{i}$--almost equal to $X_{i}$%
. Now consider the $G$--equivariant bijection $E\rightarrow E^{\prime }$
determined by sending $X_{i}$ to $Z_{i}$. The above argument shows that if $%
U $ is any element of $E$, and $U^{\prime }$ is the corresponding element of 
$E^{\prime }$, then $U\;$and $U^{\prime }$ are $stab(U)$--almost equal. We
will show that in most cases, this bijection automatically preserves the
partial orders on $E$ and $E^{\prime }$, implying that $T$ and $T^{\prime }$
are $G$--isomorphic, as required. We compare the partial orders on $E$ and $%
E^{\prime }$ rather as in the proof of Lemma \ref{equalsplittings}.

For any elements $U$ and $V$ of $E$, let $U^{\prime }$ and $V^{\prime }$
denote the corresponding elements of $E^{\prime }$. Thus $U\cap V$ is small
if and only if $U^{\prime }\cap V^{\prime }$ is small. We would like to show
that when we compare the four sets $U^{(\ast )}\cap V^{(\ast )}$ with the
four sets $U^{\prime (\ast )}\cap V^{\prime (\ast )}$, then corresponding
sets are empty, so that the partial orders are preserved by our bijection.
Otherwise, there must be $U$ and $V$ in $E$ such that two of the sets $%
U^{(\ast )}\cap V^{(\ast )}$ are small. If $U\;$and $V$ are translates of $%
X_{i}$ and $X_{j}$, then Lemma \ref{equalsplittings} tells us that the
splittings $s_{i}$ and $s_{j}$ are conjugate. As we are assuming that
distinct splittings are not conjugate, it follows that $i=j$. Now the
arguments in the proof of Lemma \ref{equalsplittings} show that either the
splitting $s_{i}$ is an amalgamated free product of the form $L\ast _{H}C$,
with $\left| L:H\right| =2$, or $H$ is normal in $G$ and $H\backslash G$ is $%
\mathbb{Z}$ or $\mathbb{Z}_{2}\ast \mathbb{Z}_{2}$. If the second case
occurs, then there can be only one splitting in the given family, so it is
immediate that $\Gamma $ and $\Gamma ^{\prime }$ have the same underlying
graph, and that corresponding edge and vertex groups are conjugate. If the
first case occurs and the partial orders on translates of $X_{i}$ and $Z_{i}$
do not match, we must have $lX_{i}^{\ast }\subset X_{i}$ but $Z_{i}\subset
lZ_{i}^{\ast }$, for all $l\in L-H$. We now pick $l\in L-H$ and alter our
bijection from $E$ to $E^{\prime }$ so that $X_{i}$ maps to $%
W_{i}=lZ_{i}^{\ast }$ and extend $G$--equivariantly to the translates of $%
X_{i}$ and $X_{i}^{\ast }$. This ensures that the partial orders on $E$ and $%
E^{\prime }$ match for translates of $X_{i}$. By repeating this for other
values of $i$ as necessary, we can arrange that the partial orders match
completely, and can then conclude that $T$ and $T^{\prime }$ are $G$%
--isomorphic as required.

We end by discussing the case when some of the given $n$ splittings are
conjugate. We divide the splittings into conjugacy classes and discard all
except one splitting from each conjugacy class, to obtain $k$ splittings.
Now we apply the preceding argument to express $G$ uniquely as the
fundamental group of a graph $\Gamma $ of groups with $k$ edges. If an edge
of $\Gamma $ corresponds to a splitting over a subgroup $H$ which is
conjugate to $r-1$ other splittings, we simply subdivide this edge into $r$
sub-edges, and label all the sub-edges and the $r-1$ new vertices by $H$.
This shows the existence of the required graph of groups structure $\Gamma
^{\prime }$ corresponding to the original $n$ splittings. The uniqueness of $%
\Gamma ^{\prime }$ follows from the uniqueness of $\Gamma $, and the fact
that the collection of all the edges of $\Gamma ^{\prime }$ which correspond
to a given splitting of $G$ must form an interval in $\Gamma ^{\prime }$ in
which all the interior vertices have valence $2$. This completes the proof
of Theorem \ref{disjointsplittings}.
\end{proof}

Now we turn to the proof of Theorem \ref{splittingsexist} that splittings
exist. It will be convenient to make the following definitions. We will use $%
H^{g}$ to denote $gHg^{-1}.$

\begin{definition}
If $X$ is a $H$--almost invariant subset of $G\;$and $Y$ is a $K$--almost
invariant subset of $G$, and if $X$ and $Y$ are $H$--almost equal, then we
will say that $X$ and $Y$ are equivalent and write $X\sim Y$. (Note that $H$
and $K$ must be commensurable.)
\end{definition}

\begin{definition}
If $H$ is a subgroup of a group $G$, the commensuriser in $G$ of $H$
consists of those elements $g$ in $G$ such that $H$ and $H^{g}$ are
commensurable subgroups of $G$. The commensuriser is clearly a subgroup of $%
G $ and is denoted by $Comm_{G}(H)$ or just $Comm(H)$, when the group $G$ is
clear from the context.
\end{definition}

Now we come to the proof of our second main result.

\begin{theorem}
\label{splittingsexist}Let $G$ be a finitely generated group with a finitely
generated subgroup $H$, such that $e(G,H)\geq 2$. If there is a non-trivial $%
H$--almost invariant subset $X$ of $G$ such that $i(H\backslash X,H\backslash
X)=0$, then $G$ has a splitting over some subgroup $H^{\prime }$
commensurable with $H$. Further, one of the $H^{\prime }$--almost invariant
sets $Y$ determined by the splitting is equivalent to $X.$
\end{theorem}

\begin{remark}
This is the best possible result of this type, as it is clear that one
cannot expect to obtain a splitting over $H$ itself. For example, suppose
that $H$ is carried by a proper power of a two-sided simple closed curve on
a closed surface whose fundamental group is $G,$ so that $e(G,H)=2.$ There
are essentially only two non-trivial almost invariant subsets of 
$H\backslash G,$ each with vanishing self-intersection number, but there is
no splitting of $G$ over $H.$
\end{remark}

\begin{proof}
The idea of the proof is much as before. We let $P$ denote the almost
invariant subset $H\backslash X$ of $H\backslash G$, and let $E$ denote $%
\{gX,gX^{\ast }:g\in G\}$. We want to apply the first part of Theorem \ref
{algebrafromTorusTheorem}. As before, the assumption that $i(P,P)=0$ implies
that $E$ is almost nested. However, in order to apply Theorem \ref
{algebrafromTorusTheorem}, we also need to know that for any pair of
elements $U$ and $V$ of $E$, if two of the four sets $U^{(\ast )}\cap
V^{(\ast )}$ are small then one is empty. In the proof of Theorem \ref
{disjointsplittings}, we simply applied Lemma \ref{emptyorconjugate}.
However, here the situation is somewhat more complicated. Lemma \ref
{Kcommensurises} below shows that if $X\cap gX^{\ast }$ and $gX\cap X^{\ast }
$ are both small, then $g$ must lie in a certain subgroup $\mathcal{K}$ of $%
Comm_{G}(H)$. Thus it would suffice to arrange that $E$ is nested with
respect to $\mathcal{K}$, ie, that $gX$ and $X$ are nested so long as $g$
lies in $\mathcal{K}$. Now Proposition \ref{almostnestedimpliesnested} below
tells us that there is a subgroup $H^{\prime }$ commensurable with $H$ and a 
$H^{\prime }$--almost invariant set $Y$ equivalent to $X$ such that $%
E^{\prime }=\{gY,gY^{\ast }:g\in G\}$ is nested with respect to $\mathcal{K}$%
. It follows that if $U\;$and $V$ are any elements of $E^{\prime }$ and if $%
U\cap V^{\ast }$ and $V\cap U^{\ast }$ are both small, then one of them is
empty. We also claim that, like $E$, the set $E^{\prime }$ is almost nested.
This means that if we let $P^{\prime }$ denote $H^{\prime }\backslash Y$, we
are claiming that $i(P^{\prime },P^{\prime })=0$. Let $H^{\prime \prime }$
denote $H\cap H^{\prime }$. The fact that $Y$ is equivalent to $X$ means
that the pre-images in $H^{\prime \prime }\backslash G$ of $P$ and of $%
P^{\prime }$ are almost equal almost invariant sets which we denote by $Q$
and $Q^{\prime }$. If $d$ denotes the index of $H^{\prime }$ in $H$, then $%
i(Q,Q)=d^{2}i(P,P)=0$ and similarly $i(Q^{\prime },Q^{\prime })$ is an
integral multiple of $i(P^{\prime },P^{\prime })$. As $Q$ and $Q^{\prime }$
are almost equal, it follows that $i(Q^{\prime },Q^{\prime })=i(Q,Q)$, and
hence that $i(P^{\prime },P^{\prime })=0$ as claimed. This now allows us to
apply Theorem \ref{algebrafromTorusTheorem} to the set $E^{\prime }$. We
conclude that $G\;$splits over the stabiliser $H^{\prime \prime }$ of $Y$,
that $H^{\prime \prime }$ contains $H^{\prime }$ with finite index and that
one of the $H^{\prime \prime }$--almost invariant sets associated to the
splitting is equivalent to $X^{\prime }$. It follows that $H^{\prime \prime }
$ is commensurable with $H$ and that one of the $H^{\prime \prime }$--almost
invariant sets determined by the splitting is equivalent to $X.$ This
completes the proof of Theorem \ref{splittingsexist} apart from the proofs
of Lemma \ref{Kcommensurises} and Proposition \ref{almostnestedimpliesnested}%
.
\end{proof}

It remains to prove the two results we just used. The proofs do not use the
hypothesis that the set of all translates of $X$ and $X^{\ast }$ are almost
nested. Thus for the rest of this section, we will consider the following
general situation.

Let $G$\ be a finitely generated group with a finitely generated subgroup $H$%
\ such that $e(G,H)\geq 2$, and let $X$ denote a non-trivial $H$--almost
invariant subset of $G.$

Recall that our problem in the proof of Theorem \ref{splittingsexist} is the
possibility that two of the four sets $X^{(\ast )}\cap gX^{(\ast )}$ are
small. As this would imply that $gX\sim X$ or $X^{\ast }$, it is clear that
the subgroup $\mathcal{K}$ of $G$ defined by $\mathcal{K}=\{g\in G:gX\sim X$
or $X^{\ast }\}$ is very relevant to our problem. We will consider this
subgroup carefully. Here is the first result we quoted in the proof of
Theorem \ref{splittingsexist}.

\begin{lemma}
\label{Kcommensurises}If $\mathcal{K}=\{g\in G:gX\sim X$ or $X^{\ast }\}$,
then $H\subset \mathcal{K}\subset Comm_{G}(H).$
\end{lemma}

\begin{proof}
The first inclusion is clear. The second is proved in essentially the same
way as the proof of the first part of Lemma \ref{equalsplittings}. Let $g$
be an element of $\mathcal{K}$, and consider the case when $gX\sim X$ ( the
other case is similar). Recall that this means that the sets $X\cap gX^{\ast
}$ and $X^{\ast }\cap gX$ are both small. Now for each edge of $\delta gX$,
either it is also an edge of $\delta X$ or it has (at least) one end in one
of the two small sets. As the images in $H\backslash \Gamma $ of $\delta X$
and of each small set is finite, and as the graph $\Gamma $ is locally
finite, it follows that the image of $\delta gX$ in $H\backslash \Gamma $
must be finite. This implies that $H\cap H^{g}$ has finite index in the
stabiliser $H^{g}$ of $\delta gX$. By reversing the roles of $X$ and $gX$,
it follows that $H\cap H^{g}$ has finite index in $H$, so that $H$ and $H^{g}
$ must be commensurable, as claimed. It follows that $\mathcal{K}\subset
Comm_{G}(H)$, as required.
\end{proof}

Another way of describing our difficulty in applying Theorem \ref
{algebrafromTorusTheorem} is to say that it is caused by the fact that the
translates of $X$ and $X^{\ast }$ may not be nested. However, Lemma \ref
{finitenumberofdoublecosets} assures us that ``most'' of the translates are
nested. The following result gives us a much stronger finiteness result.

\begin{lemma}
\label{finitelymanynonnested} Let $G$, $H$, $X$, $\mathcal{K}$ be as above.
Then $\{g\in \mathcal{K}:gX$ and $X$ are not nested$\}$ consists of a finite
number of right cosets $gH$ of $H$ in $G.$
\end{lemma}

\begin{proof}
Lemma \ref{finitenumberofdoublecosets} tells us that the given set is
contained in the union of a finite number of double cosets $HgH$. If $k\in 
\mathcal{K}$, we claim that the double coset $HkH$ is itself the union of
only finitely many cosets $gH$, which proves the required result. To prove
our claim, recall that $k^{-1}Hk$ is commensurable with $H$. Thus $k^{-1}Hk$
can be expressed as the union of cosets $g_{i}(k^{-1}Hk\cap H)$, for $1\leq
i\leq n.$ Hence 
$$
HkH=k(k^{-1}Hk)H=k\left( \cup _{i=1}^{n}g_{i}(k^{-1}Hk\cap H)\right)
H=k\left( \cup _{i=1}^{n}g_{i}H\right) =\cup _{i=1}^{n}kg_{i}H,
$$
so that $HkH$ is the union of finitely many cosets $gH$ as claimed.
\end{proof}

Now we come to the key result.

\begin{lemma}
\label{Kisunionofnormalisers}Let $G$, $H$, $X$, $\mathcal{K}$ be as above.
Then there are a finite number of finite index subgroups $H_{1},...,H_{m}$
of $H$, such that $\mathcal{K}$ is contained in the union of the groups $%
N(H_{i})$, $1\leq i\leq m$, where $N(H_{i})$ denotes the normaliser of $H_{i}
$ in $G.$
\end{lemma}

\begin{proof}
Consider an element $g$ in $\mathcal{K}$. Lemma \ref{Kcommensurises} tells
us that $H$ and $H^{g}$ are commensurable subgroups of $G$. Let $L$ denote
their intersection and let $L^{\prime }$ denote the intersection of the
conjugates of $L$ in $H$. Thus $L^{\prime }$ is of finite index in $H$ and $%
H^{g}$ and is normal in $H$. Now consider the quotient $L^{\prime
}\backslash G$. Let $P$ and $Q$ denote the images of $X$ and $gX$
respectively in $L^{\prime }\backslash G$. As before, $P$ and $Q$ are almost
invariant subsets of $L^{\prime }\backslash G$ which are almost equal or
almost complementary. Now consider the action of $L^{\prime }\backslash H$
on the left on $L^{\prime }\backslash G$. If $h$ is in $H$, then $hP=P$, so
that $hQ\overset{a}{=}Q.$ If $h(gX)$ and $gX$ are nested, there are four
possible inclusions, but the fact that $hQ\overset{a}{=}Q$ excludes two of
them. Thus we must have $hQ\subset Q$ or $Q$ $\subset hQ$. This implies that 
$hQ=Q$ as some power of $h$ lies in $L^{\prime }$ and so acts trivially on $%
L^{\prime }\backslash G.$ We conclude that if $h$ is an element of $%
H-L^{\prime }$ such that $h(gX)$ and $gX$ are nested, then $h$ stabilises $gX
$ and so lies in $H^{g}$. Hence $h$ lies in $L$. It follows that for each
element $h$ of $H-L$, the sets $h(gX)$ and $gX$ are not nested. Recall from
Lemma \ref{finitelymanynonnested} that $\{g\in \mathcal{K}:gX$ and $X$ are
not nested$\}$ consists of a finite number of cosets $gH$ of $H$ in $G$. It
will be convenient to denote this number by $d-1$. Thus, for $g\in \mathcal{K%
}$, the set $\{h\in \mathcal{K}:h(gX)$ and $gX$ are not nested$\}$ consists
of $d-1$ cosets $hH^{g}$ of $H^{g}$ in $G$. It follows that $H-L$ lies in
the union of $d-1$ cosets $hH^{g}$ of $H^{g}$ in $G$. As $L=H\cap H^{g}$, it
follows that $H-L$ lies in the union of $d-1$ cosets $hL$ of $L$ in $G$ and
hence that $L$ has index at most $d$ in $H.$

A similar argument shows also that $L$ has index at most $d$ in $H^{g}$. Of
course, the same bound applies to the index of $H\cap H^{g^{i}}$ in $H$, for
each $i$. Now we define $H^{\prime }=\cap _{i\in \mathbb{Z}}H^{g^{i}}.$
Clearly $H^{\prime }$ is a subgroup of $H$ which is normalised by $g$. Now
each intersection $H\cap H^{g^{i}}$ has index at most $d$ in $H$, and so $%
H^{\prime }=\cap _{i\in \mathbb{Z}}\left( H\cap H^{g^{i}}\right) $ is an
intersection of subgroups of $H$ of index at most $d$. If $H$ has $n$
subgroups of index at most $d$, it follows that $H^{\prime }$ has index at
most $d^{n}$ in $H$. Hence each element of $\mathcal{K}$ normalises a
subgroup of $H$ of index at most $d^{n}$ in $H$. As $H$ has only finitely
many such subgroups, we have proved that there are a finite number of finite
index subgroups $H_{1},...,H_{m}$ of $H$, such that $\mathcal{K}$ is
contained in the union of the groups $N(H_{i})$, $1\leq i\leq m$, as
required.
\end{proof}

Using this result, we can prove the following.

\begin{lemma}
Let $G$, $H$, $X$, $\mathcal{K}$ be as above. Then there is a subgroup $%
H^{\prime }$ of finite index in $H$, such that $\mathcal{K}$ normalises $%
H^{\prime }.$
\end{lemma}

\begin{proof}
We will consider how $\mathcal{K}$ can intersect the normaliser of a
subgroup of finite index in $H$. Let $H_{1}$ denote a subgroup of $H$ of
finite index. We denote the image of $X$ in $H_{1}\backslash G$ by $P$. Then 
$P$ is an almost invariant subset of $H_{1}\backslash G$. We consider the
group $\mathcal{K}\cap N(H_{1})$, which we will denote by $\mathcal{K}_{1}$.
Then $H_{1}\backslash \mathcal{K}_{1}$ acts on the left on $H_{1}\backslash G
$, and we have $kP\overset{a}{=}P$ or $P^{\ast }$, for every element $k$ of $%
H_{1}\backslash \mathcal{K}_{1}$, because every element of $\mathcal{K}$
satisfies $kX\sim X$ or $X^{\ast }$. Now we apply Theorem 5.8 from \cite
{Scott-Wall:Topological} to the action of $H_{1}\backslash \mathcal{K}_{1}$
on the left on the graph $H_{1}\backslash \Gamma $. This result tells us
that if $H_{1}\backslash \mathcal{K}_{1}$ is infinite, then it has an
infinite cyclic subgroup of finite index. Further the proof of this result
in \cite{Scott-Wall:Topological} shows that the quotient of $H_{1}\backslash
\Gamma $ by $H_{1}\backslash \mathcal{K}_{1}$ must be finite. This implies
that $H_{1}\backslash \Gamma $ has two ends and that $\mathcal{K}_{1}$ has
finite index in $G$. Hence either $H_{1}\backslash \mathcal{K}_{1}$ is
finite, or it has two ends and $\mathcal{K}_{1}$ has finite index in $G.$

Recall that there are a finite number of finite index subgroups $%
H_{1},...,H_{m}$ of $H$, such that $\mathcal{K}$ is contained in the union
of the groups $N(H_{i})$, $1\leq i\leq m$. The above discussion shows that,
for each $i$, if $\mathcal{K}_{i}$ denotes $\mathcal{K}\cap N(H_{i})$,
either $H_{i}\backslash \mathcal{K}_{i}$ is finite, or it has two ends and $%
\mathcal{K}_{i}$ has finite index in $G$. We consider two cases depending on
whether or not every $H_{i}\backslash \mathcal{K}_{i}$ is finite.

Suppose first that each $H_{i}\backslash \mathcal{K}_{i}$ is finite. We
claim that $\mathcal{K}$ contains $H$ with finite index. To see this, let $%
H^{\prime \prime }=\cap H_{i}$, so that $H^{\prime \prime }$ is a subgroup
of $H$ of finite index, and note that $\mathcal{K}$ is the union of a finite
collection of groups $\mathcal{K}_{i}$ each of which contains $H^{\prime
\prime }$ with finite index, so that $\mathcal{K}$ is the union of finitely
many cosets of $H^{\prime \prime }$. It follows that $\mathcal{K}$ also
contains $H^{\prime \prime }$ with finite index and hence contains $H$ with
finite index as claimed. If we let $H^{\prime }$ denote the intersection of
the conjugates of $H$ in $\mathcal{K}$, then $H^{\prime }$ is the required
subgroup of $H$ which is normalised by $\mathcal{K}.$

Now we turn to the case when $H_{1}\backslash \mathcal{K}_{1}$ is infinite
and so $H_{1}\backslash \mathcal{K}_{1}$ has two ends and $\mathcal{K}_{1}$
has finite index in $G.$ Define $H^{\prime }$ to be $\cap _{k\in \mathcal{K}%
}\left( H_{1}\right) ^{k}$. As $\mathcal{K}$ contains $\mathcal{K}_{1}$ with
finite index, $H^{\prime }$ is the intersection of only finitely many
conjugates of $H_{1}$. As $\mathcal{K}$ is contained in $Comm(H)$, each of
these conjugates of $H_{1}$ is commensurable with $H_{1}$. It follows that $%
H^{\prime }$ is a subgroup of $H$ of finite index in $H$ which is normalised
by $\mathcal{K}$. This completes the proof of the lemma.
\end{proof}

The key point here is that $\mathcal{K}$ normalises $H^{\prime }$ rather
than just commensurises it. Now we can prove the second result which we
quoted in the proof of Theorem \ref{splittingsexist}.

\begin{proposition}
\label{almostnestedimpliesnested}Suppose that $(G,H)$ is a pair of finitely
generated groups and that $X$ is a non-trivial $H$--almost invariant subset
of $G$. Then, there is a subgroup $H^{\prime }$ of $G$ which is
commensurable with $H$, and a non-trivial $H^{\prime }$--almost invariant set 
$Y$ equivalent to $X$ such that $\{gY,gY^{\ast }:g\in G\}$ is nested with
respect to the subgroup $\mathcal{K}=\{g\in G:gX\sim X$ or $X^{\ast }\}$ of $%
G$.
\end{proposition}

\begin{proof}
The previous lemma tells us that there is a subgroup $H^{\prime }$ of finite
index in $H$ such that $\mathcal{K}$ normalises $H^{\prime }$. Let $P$
denote the almost invariant subset $H\backslash X$ of $H\backslash G$, and
let $P^{\prime }$ denote the almost invariant subset $H^{\prime }\backslash X
$ of $H^{\prime }\backslash G.$

Suppose that the index of $H^{\prime }$ in $\mathcal{K}$ is infinite. Recall
from the proof of the preceding lemma that $H^{\prime }\backslash \mathcal{K}
$ has two ends and that $\mathcal{K}$ has finite index in $G$. We construct
a new non-trivial $H^{\prime \prime }$--almost invariant set $Y$ as follows.
Since the quotient group $H^{\prime }\backslash \mathcal{K}$ has two ends, $%
\mathcal{K}$ splits over a subgroup $H^{\prime \prime }$ which contains $%
H^{\prime }$ with finite index. Thus there is a $H^{^{\prime \prime }}$%
--almost invariant set $X^{\prime \prime }$ in $\mathcal{K}$ which is nested
with respect to $\mathcal{K}$. Further, $H^{\prime \prime }$ is normal in $%
\mathcal{K}$ and the quotient group must be isomorphic to $\mathbb{Z}$ or $%
\mathbb{Z}_{2}\ast \mathbb{Z}_{2}$. Let $\{g_{1}=e,g_{2},...,g_{n}\}$ be
coset representatives of $\mathcal{K}$ in $G$ so that $G=\cup _{i}\mathcal{K}%
g_{i}$. We take $Y=\cup _{i}X^{\prime \prime }g_{i}$. It is easy to check
that $Y$ is $H^{\prime \prime }$--almost invariant and that $\{gY,gY^{\ast
}:g\in G\}$ is nested with respect to $\mathcal{K}.$

Now suppose that the index of $H^{\prime }$ in $\mathcal{K}$ is finite. We
will define the subgroup $\mathcal{K}_{0}=\{g\in G:gX\sim X\}$ of $\mathcal{K%
}$. The index of $\mathcal{K}_{0}$ in $\mathcal{K}$ is at most two.

First we consider the case when $\mathcal{K}=\mathcal{K}_{0}$. We define $%
P^{\prime \prime }$ to be the intersection of the translates of $P^{\prime }$
under the action of $H^{\prime }\backslash \mathcal{K}$. Thus $P^{\prime
\prime }$ is invariant under the action of $H^{\prime }\backslash \mathcal{K}
$. As all the translates of $P^{\prime }$ by elements of $H^{\prime
}\backslash \mathcal{K}$ are almost equal to $P^{\prime }$, it follows that $%
P^{\prime \prime }\overset{a}{=}P^{\prime }$ so that $P^{\prime \prime }$ is
also an almost invariant subset of $H^{\prime }\backslash G$. Let $Y$ denote
the inverse image of $P^{\prime \prime }$ in $G$, so that $Y$ is invariant
under the action of $\mathcal{K}$. In particular, $\{gY,gY^{\ast }:g\in G\}$
is nested with respect to $\mathcal{K}$, as required.

Now we consider the general case when $\mathcal{K}\neq \mathcal{K}_{0}$. We
can apply the above arguments using $\mathcal{K}_{0}$ in place of $\mathcal{K%
}$ to obtain a subgroup $H^{\prime \prime }$ of $G$ and a $H^{\prime \prime }
$--almost invariant subset $Y$ of $G$ which is equivalent to $X$, and whose
translates are nested with respect to $\mathcal{K}_{0}$. We also know that $Y
$ is $\mathcal{K}_{0}$--invariant. Let $Q$ denote the image of $Y$ in $%
\mathcal{K}_{0}\backslash G$, let $k$ denote an element of $\mathcal{K}-%
\mathcal{K}_{0}$ and consider the involution of $\mathcal{K}_{0}\backslash G$
induced by $k$. Then $Q$ is a non-trivial almost invariant subset of $%
\mathcal{K}_{0}\backslash G$ and $kQ\overset{a}{=}Q^{\ast }$. Define $R=Q-kQ$%
, so that $R\overset{a}{=}Q$ and let $Z$ denote the pre-image of $R$ in $G$.
We claim that the translates of $Z$ and $Z^{\ast }$ are nested with respect
to $\mathcal{K}$. First we show that they are nested with respect to $%
\mathcal{K}_{0}$, by showing that $Z=Y-kY$ is $\mathcal{K}_{0}$--invariant.
For $k_{0}\in \mathcal{K}_{0}$, we have $k^{-1}k_{0}k\in \mathcal{K}_{0}$ as 
$\mathcal{K}_{0}$ must be normal in $\mathcal{K}$. It follows that $%
k_{0}kY=kY$. As $k_{0}Y=Y$, we see that $Z$ is $\mathcal{K}_{0}$--invariant
as required. In order to show that the translates of $Z$ and $Z^{\ast }$ are
nested with respect to $\mathcal{K}$, we will also show that $Z\cap kZ$ is
empty. This follows from the fact that $R\cap kR=\left( Q-kQ\right) \cap
k\left( Q-kQ\right) =\left( Q-kQ\right) \cap \left( kQ-Q\right) $ which is
clearly empty.

This completes the proof of Proposition \ref{almostnestedimpliesnested}.
\end{proof}

\section{Strong intersection numbers}

Let $G$ be a finitely generated group and let $H$ and $K$ be subgroups of $G$%
. Let $X$ be a non-trivial $H$--almost invariant subset of $G$ and let $Y$ be
a non-trivial $K$--almost invariant subset of $G.$ In section 1, we discussed
what it means for $X$ to cross $Y$ and the fact that this is symmetric. As
mentioned in the introduction, there is an alternative way to define
crossing of almost invariant sets. Recall that, in section 1, we introduced
our definition of crossing by discussing curves on surfaces. Thus it seems
natural to discuss the crossing of $X$ and $Y$ in terms of their boundaries.
We call this strong crossing. However, this leads to an asymmetric
intersection number. In this section, we define strong crossing and discuss
its properties and some applications.

We consider the Cayley graph $\Gamma $ of $G$ with respect to a finite
system of generators. We will usually assume that $H$ and $K$ are finitely
generated though this does not seem necessary for most of the definitions
below. We will also think of $\delta X$ as a set of edges in $\Gamma $ or as
a set of points in $G$, where the set of points will simply be the
collection of endpoints of all the edges of $\delta X.$

\begin{definition}
We say that $Y$\textsl{\ crosses }$X$\textsl{\ strongly} if both $\delta
Y\cap X$ and $\delta Y\cap X^{\ast }$ project to infinite sets in $%
H\backslash G$.
\end{definition}

\begin{remark}
\label{strongcrossingimpliescrossing}This definition is independent of the
choice of generators for $G$ which is used to define $\Gamma $. Clearly, if $%
Y$ crosses $X$ strongly, then $Y$ crosses $X$.
\end{remark}

Strong crossing is not symmetric. For an example, one need only consider an
essential two-sided simple closed curve $S$ on a compact surface $F$ which
intersects a simple arc $L$ transversely in a single point. Let $G$ denote $%
\pi _{1}(F)$, and let $H$ and $K$ respectively denote the subgroups of $G$
carried by $S$ and $L$, so that $H$ is infinite cyclic and $K$ is trivial.
Then $S$ and $L$ each define a splitting of $G$ over $H\;$and $K$
respectively. Let $X$ and $Y$ denote associated standard $H$--almost
invariant and $K$--almost invariant subsets of $G$. These correspond to
submanifolds of the universal cover of $F$ bounded respectively by a line $%
\widetilde{S}$ lying above $S$ and by a compact interval $\widetilde{L}$
lying above $L$, such that $\widetilde{S}$ meets $\widetilde{L}$
transversely in a single point. Clearly, $X$ crosses $Y$ strongly but $Y$
does not cross $X$ strongly.

However, a strong intersection number can be defined as before. It is
usually asymmetric, but we will be particularly interested in the case of
self-intersection numbers when this asymmetry will not arise.

\begin{definition}
The strong intersection number $si(H\backslash X,K\backslash Y)$ is defined
to be the number of double cosets $KgH$ such that $gX$ crosses $Y$ strongly.
In particular, $si(H\backslash X,H\backslash X)=0$ if and only if at least
one of $\delta gX\cap X$ and $\delta gX\cap X^{\ast }$ is $H$--finite, for
each $g\in G$.
\end{definition}

\begin{remark}
If $s$ and $t$ are splittings of a group $G$ over subgroups $H$ and $K$,
with associated almost invariant subsets $X$ and $Y$ of $G$, it is natural
to say that $s$ crosses $t$ strongly if $si(H\backslash X,K\backslash Y)\neq
0$. It is easy to show that this is equivalent to the idea introduced by
Sela \cite{Sela:JSJ} that $s$ is hyperbolic with respect to $t$.
\end{remark}

Remark \ref{strongcrossingimpliescrossing} shows that $si(H\backslash
X,H\backslash X)\leq i(H\backslash X,H\backslash X)$. Recall that Theorem 
\ref{splittingsexist} shows that if $i(H\backslash X,H\backslash X)=0$, then 
$G$ splits over a subgroup $H^{\prime }$ commensurable with $H$. Thus the
vanishing of the strong self-intersection number may be considered as a
first obstruction to splitting $G$ over some subgroup related to $H$. We
will show in Corollary \ref{canmakenested} that the vanishing of the strong
self-intersection number has a nice algebraic formulation. This is that when 
$si(H\backslash X,H\backslash X)$ vanishes, we can find a subgroup $K$ of $G$%
, commensurable with $H$, and a $K$--almost invariant subset $Y$ of $G$ which
is nested with respect to $Comm_{G}(H)=Comm_{G}(K)$. However, $Y$ may be
very different from $X$. This leads to some splitting results when we place
further restrictions on $H$.

\begin{proposition}
\label{conditionforsi=0}Let $G$ be a finitely generated group with finitely
generated subgroup $H$, and let $X$ be a non-trivial $H$--almost invariant
subset of $G$. Then $si(H\backslash X,H\backslash X)=0$ if and only if there
is a subset $Y$ of $G$ which is $H$--almost equal to $X$ (and hence $H$%
--almost invariant) such that $HYH=Y$.
\end{proposition}

\begin{proof}
Suppose that there exists a subset $Y$ of $G$ which is $H$--almost equal to $%
X $, such that $HYH=Y$. We have 
\begin{equation*}
si(H\backslash X,H\backslash X)=si(H\backslash Y,H\backslash Y),
\end{equation*}
as $X$ and $Y$ are $H$--almost equal. So, it is enough to show that for every 
$g\in G$, either $g\delta Y\cap Y$ or $g\delta Y\cap Y^{\ast }$ is $H$%
--finite. Suppose that $g\in Y$. Consider $\delta Y\cap Y$ which is a union
of a finite number of right cosets $Hg_{i}$, $1\leq i\leq n$. Since $g\in Y$%
, $gH\subset Y$. For any $h\in H$, $d(gh,ghg_{i})=d(1,g_{i})$. Thus $g\delta
Y$ is at a bounded distance from $Y$ and hence $g\delta Y\cap Y^{\ast }$ has
finite image in $H\backslash G$. Similarly, if $g\in Y^{\ast }$, $g\delta
Y\cap Y$ projects to a finite set in $H\backslash G$.

For the converse, suppose that $si(H\backslash X,H\backslash X)=0$ and let $%
\pi $ denote the projection from $G$ to $H\backslash G$. By hypothesis, $\pi
(g\delta X)\cap (H\backslash X)$ or $\pi (g\delta X)\cap (H\backslash
X^{\ast })$ is finite. The proof of Lemma \ref{posatisfiesDunwoody} tells us
that there is a positive number $d$ such that, for every $g\in G$, the set $%
g\delta X$ is contained in a $d$--neighbourhood of $X$ or $X^{\ast }$. Let $%
V=N(X,d)$, the $d$--neighbourhood of $X$ and let $Y=\{g|g(\delta X)\subset
V\} $. If $g\in Y$ and $h\in H$, then $hg\delta X\subset hV=V$ and thus $%
HY=Y $. If $g\in Y$ and $h\in H$, then $gh(\delta X)=g(\delta X)\subset V$
and thus $YH=Y$. It only remains to show that $Y$ is $H$--almost equal to $X.$
This is essentially shown in the third and fourth paragraphs of the proof of
Theorem \ref{algebrafromTorusTheorem}.
\end{proof}

\begin{definition}
We will say that a pair of finitely generated groups $(G,H)$ is of \emph{%
surface type} if $e(G,H^{\prime })=2$ for every subgroup $H^{\prime }$ of
finite index in $H$ and $e(G,H^{\prime })=1$ for every subgroup $H^{\prime }$
of infinite index in $H$.
\end{definition}

This terminology is suggested by the dichotomy in \cite{S-S}. Note that for
such pairs any two non-trivial $H$--almost invariant sets in $G$ are $H$%
--almost equal or $H$--almost complementary. We will see that for pairs of
surface type, strong and ordinary intersection numbers are equal.

\begin{proposition}
\label{crossingiffstrongcrossing}Let $(G,H)$ be a pair of surface type, let $%
X$ be a non-trivial $H$--almost invariant subset of $G$ and let $Y$ be a
non-trivial $K$--almost invariant subset of $G$ for some subgroup $K$ of $G$.
Then $Y$ crosses $X$ if and only if $Y$ crosses $X$ strongly.
\end{proposition}

\begin{proof}
Let $\Gamma $ be the Cayley graph of $G$ with respect to a finite system of
generators and let $P=H\backslash X$. As in the proof of Lemma \ref
{finitenumberofdoublecosets}, for a set $S$ of vertices in a graph, we let $%
\overline{S}$ denote the maximal subgraph with vertex set equal to $S$. We
will show that exactly one component of $\overline{X}$ has infinite image in 
$H\backslash \Gamma $. Note that $\overline{P}$ has exactly one infinite
component as $H\backslash \Gamma $ has only two ends. Let $Q$ denote the set
of vertices of the infinite component of $\overline{P}$ and let $W$ denote
the inverse image of $Q$ in $G$. If $\overline{W}$ has components with
vertex set $L_{i}$, then we have $\cup \delta (L_{i})=\delta W\subseteq
\delta X$. Let $L$ denote the vertex set of a component of $\overline{W}$,
and let $H_{L}$ be the stabilizer in $H$ of $L$. Since $\delta Q$ is finite,
we see that $H_{L}\backslash \delta L$ is finite. Hence $H_{L}\backslash
\Gamma $ has more than one end. Now our hypothesis that $(G,H)$ is of
surface type implies that $H_{L}$ has finite index in $H$ and thus $%
H_{L}\backslash \delta W$ is finite. If $H_{L}\neq H$, we see that $%
H_{L}\backslash \delta W$ divides $H_{L}\backslash \Gamma $ into at least
three infinite components. Thus $H_{L}=H$ and so $\overline{W}$ is
connected. The other components of $\overline{X}$ have finite image in $%
H\backslash \Gamma $. Similarly, exactly one component of $\overline{X^{\ast
}}$ has infinite image in $H\backslash \Gamma $. The same argument shows
that for any finite subset $D$ of $H\backslash \Gamma $ containing $\delta P$%
, the two infinite components of $((H\backslash \Gamma )-D)\cap P$ and$%
((H\backslash \Gamma )-D)\cap P^{\ast }$ have connected inverse images in $%
\Gamma $.

Recall that if $Y$ crosses $X$ strongly, then $Y$ crosses $X$. We will next
show that if $Y$ does not cross $X$ strongly, then $Y$ does not cross $X$.
Suppose that $\delta Y\cap X$ projects to a finite set in $H\backslash
\Gamma $. Take a compact set $D$ in $H\backslash \Gamma $ large enough to
contain $\delta Y\cap X$ and $\delta P$. By the argument above, if $R$ is
the infinite component of $((H\backslash \Gamma )-D)\cap P$, then its
inverse image $Z$ is connected and is contained in $\overline{X}$. Any two
points in $Z$ can be connected by a path in $Z$ and thus the path does not
intersect $\delta Y$. Thus $Z$ is contained in $Y$ or $Y^{\ast }$. Hence $%
Z\cap Y$ or $Z\cap Y^{\ast }$ is empty. Suppose that $Z\cap Y$ is empty.
Then $Z^{\ast }\supseteq Y$. Since $Z^{\ast }\cap X$ projects to a finite
set, we see that $Y\cap X$ projects to a finite set. Similarly, if $Z\cap
Y^{\ast }$ is empty, then $Y^{\ast }\cap X$ projects to a finite set in $%
H\backslash G$. Thus, we have shown that if $\delta Y\cap X$ projects to a
finite set, then either $Y\cap X$ or $Y^{\ast }\cap X$ projects to finite
set. Thus $Y$ does not cross $X$.
\end{proof}

From the above proposition and the fact that ordinary crossing is symmetric,
we deduce:

\begin{corollary}
If $(G,H)$ and $(G,K)$ are both of surface type and $X$ is a non-trivial $H$%
--almost invariant set in $G$, and $Y$ is a non-trivial $K$--almost invariant
set in $G$ then $si(H\backslash X,K\backslash Y)=i(H\backslash X,K\backslash
Y)$. In particular $i(H\backslash X,H\backslash X)=0$ if and only if $%
si(H\backslash X,H\backslash X)=0$.
\end{corollary}

Let $K$ be a Poincar\'{e} duality group of dimension $(n-1)$ which is a
subgroup of a Poincar\'{e} duality group $G$ of dimension $n$. Thus the pair 
$(G,K)$ is of surface type. In \cite{KR:1}, Kropholler and Roller defined an
obstruction $sing(K)$ to splitting $G$ over a subgroup commensurable with $K$%
. Their main result was that $sing(K)$ vanishes if and only if $G$ splits
over a subgroup commensurable with $K$. At an early stage in their proof,
they showed that $sing(K)$ vanishes if and only if there is a $K$--almost
invariant subset $Y$ of $G$ such that $KYK=Y$. Starting from this point,
Proposition \ref{conditionforsi=0}, the above Corollary and then Theorem \ref
{splittingsexist} give an alternative proof of their splitting result. Thus
Theorem \ref{splittingsexist} may be considered as a generalization of their
splitting theorem. We next reformulate in our language a conjecture of
Kropholler and Roller \cite{KR:2}:

\begin{conjecture}
\label{KRconjecture}If $G$ is a finitely generated group with a finitely
generated subgroup $H$, and if $X$ is a non-trivial $H$--almost invariant
subset of $G$ such that $si(H\backslash X,H\backslash X)=0$, then $G$ splits
over a subgroup commensurable with a subgroup of $H$.
\end{conjecture}

Note that Theorem \ref{splittingsexist} has a stronger hypothesis than this
conjecture, namely the vanishing of the self-intersection number $%
i(H\backslash X,H\backslash X)$, rather than the vanishing of the strong
self-intersection number, and it has a correspondingly stronger conclusion,
namely that $G$ splits over a subgroup commensurable with $H$ itself. A key
difference between the two statements is that, in the above conjecture, one
does not expect the almost invariant set associated to the splitting of $G$
to be at all closely related to $X$. Dunwoody and Roller proved this
conjecture when $H$ is virtually polycyclic \cite{D-R:poly}, and Sageev \cite
{Sag} proved it for quasiconvex subgroups of hyperbolic groups. The paper of
Dunwoody and Roller \cite{D-R:poly} contains information useful in the
general case. The second step in their proof, which uses a theorem of
Bergman \cite{Bergman}, proves the following result, stated in our language.
(There is an exposition of Bergman's argument and parts of \cite{D-R:poly}
in the later versions of \cite{D-S:torus}.)

\begin{theorem}
Let $(G,H)$ be a pair of finitely generated groups, and let $X$ be a $H$%
--almost invariant subset of $G$. If $si(H\backslash X,H\backslash X)=0$,
then there is a subgroup $H^{\prime }$ commensurable with $H$, and a
non-trivial $H^{\prime }$--almost invariant set $Y$ with $si(H^{\prime
}\backslash Y,H^{\prime }\backslash Y)=0$ such that the set $\{gY,gY^{\ast
}:g\in G\}$ is almost nested with respect to $Comm_{G}(H)=Comm_{G}(H^{\prime
})$.
\end{theorem}

This combined with Proposition \ref{almostnestedimpliesnested} gives:

\begin{corollary}
\label{canmakenested}With the hypotheses of the above theorem we can choose $%
H^{\prime }$ and a non-trivial $H^{\prime }$--almost invariant set $Y$ with $%
si(H^{\prime }\backslash Y,H^{\prime }\backslash Y)=0$ such that $%
\{gY,gY^{\ast }:g\in G\}$ is almost nested with respect to $Comm_{G}(H)$ and
is nested with respect to the subgroup $\mathcal{K}=\{g\in G:gX\sim X$ or $%
X^{\ast }\}$ of $Comm_{G}(H)$.
\end{corollary}

Now Theorem \ref{algebrafromTorusTheorem} yields the following
generalization of Stallings' Theorem \cite{Stallings:ends} already noted by
Dunwoody and Roller \cite{D-R:poly}:

\begin{theorem}
If $G$, $H$ are finitely generated groups with $e(G,H)>1$ and if $G$
commensurises $H$, then $G$ splits over a subgroup commensurable with $H$.
\end{theorem}

Corollary \ref{canmakenested} leads to the following partial solution of the
above conjecture of Kropholler and Roller:

\begin{theorem}
If $G$, $H$ are finitely generated groups with $e(G,H)>1$, if $e(G,K)=1$ for
every subgroup $K$ commensurable with a subgroup of infinite index in $H$,
and if $X$ is a $H$--almost invariant subset of $G$ such that $si(H\backslash
X,H\backslash X)=0$, then $G$ splits over a subgroup commensurable with $H$.
\end{theorem}

\begin{proof}
Observe that Corollary \ref{canmakenested} shows that, by changing $H$ up to
commensurability, and changing $X$, we may assume that the translates of $X$
are almost nested with respect to $Comm_{G}(H)$ and nested with respect to $%
\mathcal{K}=\{g\in G:gX\sim X$ or $X^{\ast }\}$. If we do not have almost
nesting for all translates of $X$, then there is $g$ outside $Comm_{G}(H)$
such that none of $X^{(\ast )}\cap gX^{(\ast )}$ is $H$--finite. In
particular, none of these sets is $(H\cap H^{g})$--finite. But these four
sets are each invariant under $H\cap H^{g}$ and the fact that the strong
intersection number vanishes shows that at least one of them has boundary
which is $(H\cap H^{g})$--finite. Since $g$ is not in $Comm_{G}(H)$, we have
a contradiction to our hypothesis that $e(G,K)=1$ with $K=H\cap H^{g}$. This
completes the proof.
\end{proof}

We note another application of groups of surface type which provides an
approach to the Algebraic Torus Theorem \cite{D-S:torus} similar to ours in 
\cite{S-S}. We will omit a complete discussion of this approach, but will
prove the following proposition to illustrate the ideas.

\begin{proposition}
\label{bigcommensuriserimpliesvirtuallynormal}If $(G,H)$ is of surface type
and if $H$ has infinite index in $Comm_{G}(H)$, then there is a subgroup $%
H^{\prime }$ of finite index in $H$ such that the normalizer $N(H^{\prime })$
of $H^{\prime }$ is of finite index in $G$ and $H^{\prime }\backslash
N(H^{\prime })$ is virtually infinite cyclic. In particular, if $H$ is
virtually polycyclic, then $G$ is virtually polycyclic.
\end{proposition}

\begin{proof}
Let $X$ be a non-trivial $H$--almost invariant subset of $G$, let $g$ be an
element of $Comm_{G}(H)$ and let $Y=gX$, so that $Y\;$has stabiliser $H^{g}$%
. Let $H^{\prime }$ denote the intersection $H\cap H^{g}$ which has finite
index in both $H$ and in $H^{g}$ because $g$ lies in $Comm_{G}(H)$. Thus $%
H^{\prime }\backslash X$ and $H^{\prime }\backslash Y$ are both almost
invariant subsets of $H^{\prime }\backslash G$. As $(G,H)$ is of surface
type, the pair $(G,H^{\prime })$ has two ends so that $H^{\prime }\backslash
X$ and $H^{\prime }\backslash Y$ are almost equal or almost complementary.
It follows that $X$ is $H$--almost equal to $Y$ or $Y^{\ast }$, ie, $gX\sim X
$ or $gX\sim X^{\ast }$. Recall from Lemma \ref{Kcommensurises}, that if $%
\mathcal{K}$ denotes $\{g\in G:gX\sim X$ or $gX\sim X^{\ast }\}$, then $%
\mathcal{K}\subset Comm_{G}(H)$. It follows that in our present situation $%
\mathcal{K}$ must equal $Comm_{G}(H)$. By Lemma \ref{Kisunionofnormalisers},
we see that there are a finite number of subgroups $H_{1},...,H_{m}$ of
finite index in $H$ such that $\mathcal{K}$ is contained in the union of the
normalizers $N(H_{i})$. As $H$ has infinite index in $\mathcal{K}=Comm_{G}(H)
$, one of the $H_{i}$, say $H_{1}$, has infinite index in its normalizer $%
N(H_{1})$. As $(G,H)$ is of surface type, the pair $(G,H_{1})$ has two ends,
so we can apply Theorem 5.8 from \cite{Scott-Wall:Topological} to the action
of $H_{1}\backslash N(H_{1})$ on the left on the graph $H_{1}\backslash
\Gamma $. This result tells us that $H_{1}\backslash N(H_{1})$ is virtually
infinite cyclic. Further the proof of this result in \cite
{Scott-Wall:Topological} shows that the quotient of $H_{1}\backslash \Gamma $
by $H_{1}\backslash N(H_{1})$ must be finite so that $N(H_{1})$ has finite
index in $G$.
\end{proof}

The arguments of \cite{S-S} can be extended to show:

\begin{theorem}
\label{toomanyends} Let $(G,H)$ be a pair of finitely generated groups with $%
H$ virtually polycyclic and suppose that $G$ does not split over a subgroup
commensurable with a subgroup of infinite index in $H$. If for some subgroup 
$K$ of $H$, $e(G,K)\geq 3$, then $G$ splits over a subgroup commensurable
with $H$.
\end{theorem}

We end this section with an interpretation of intersection numbers in the
case when the strong and ordinary intersection numbers are equal. This
corrects a mistake in \cite{Scott:Intersectionnumbers}. Suppose that a group 
$G$ splits over subgroups $H$ and $K$ and let the corresponding $H$--almost
and $K$--almost invariant subsets of $G$ be $X$ and $Y$. Let $T$ denote the
Bass--Serre tree corresponding to the splitting of $G$ over $K$ and consider
the action of $H$ on $T$. Let $T^{\prime }$ denote the minimal $H$--invariant
subtree of $T$, and let $\Psi $ denote the quotient graph $H\backslash
T^{\prime }.$ Similarly, we get a graph $\Phi $ by considering the action of 
$K$ on the Bass--Serre tree corresponding to the splitting of $G$ over $H$.
We have:

\begin{theorem}
With the above notation, suppose that\newline
$i(H\backslash X,K\backslash Y)=si(H\backslash X,K\backslash Y)$. Then the
number of edges in $\Psi $ is the same as the number of edges in $\Phi $ and
both are equal to $si(H\backslash X,K\backslash Y)$.
\end{theorem}

\begin{proof}
The proof of Theorem 3.1 of \cite{Scott:Intersectionnumbers} goes through
because of our assumption that $i(H\backslash X,K\backslash
Y)=si(H\backslash X,K\backslash Y)$. The mistake in \cite
{Scott:Intersectionnumbers} occurs in the proof of Lemma 3.6 of \cite
{Scott:Intersectionnumbers} where it is implicitly assumed that if $X$
crosses $Y$, then it crosses $Y$ strongly. Since we have assumed that the
two intersection numbers are equal, the argument is now valid.
\end{proof}

\rk{Acknowledgement}The first author is partially supported by NSF
grants DMS 034681 and 9626537.

\end{document}